\documentclass[11pt]{elsarticle}

\usepackage{amsmath,amsthm,bm}

 \newtheorem{theorem}{Theorem}

\makeatletter
\def\ps@pprintTitle{%
 \let\@oddhead\@empty
 \let\@evenhead\@empty
 \def\@oddfoot{}%
 \let\@evenfoot\@oddfoot}
\makeatother

\usepackage{tabularx} 
\usepackage{caption}  
\usepackage[export]{adjustbox}

\usepackage{placeins,subcaption,tabularx}
\usepackage{amsmath,amscd,wrapfig,amsfonts,amssymb,mathtools,mathrsfs}
\usepackage{booktabs,multirow,lscape,cool,bm}
\usepackage{url,parskip,caption}

\usepackage[top=1.10in, bottom=1.10in, left=1.10in, right=1.10in]{geometry}

\makeatletter
\g@addto@macro\normalsize{%
  \setlength\abovedisplayskip{.6em}
  \setlength\belowdisplayskip{.6em}
  \setlength\abovedisplayshortskip{.6em}
  \setlength\belowdisplayshortskip{.6em}
}

\newcommand{\qmin}{q_{\mbox{\tiny min}}}

\newcommand{\Dk}{\mathscr{D}_k}
\newcommand{\rlg}{r_{\mbox{\tiny LG}}}

\newcommand{\tcheb}[1]{t_{#1,k}^{\mbox{\tiny cheb}}}
\newcommand{\rnon}{r_{\mbox{\tiny non}}}

\begin{document}

\begin{frontmatter}


\title
{
An accelerated frequency-independent solver for oscillatory differential equations
}

\begin{abstract}
Oscillatory second order linear ordinary differential equations 
arise in many  numerical and scientific calculations.  Because the running 
times  of standard ODE solvers increase roughly linearly with frequency when
they are applied to such problems, a variety of specialized methods,
most of them quite complicated, have been proposed.
Here, we point out that one of the simplest conceivable approaches not only
works, but yields a scheme for  solving oscillatory second order linear
ordinary differential equations which is significantly faster than 
current state-of-the-art techniques.
Our method, which operates by constructing a slowly varying phase function representing 
a basis in the space of solutions of the differential equation, runs in time independent
of the frequency of oscillations of the solutions and can be applied
to second order equations whose solutions are 
oscillatory in some regions and slowly varying in others.
In the high-frequency regime, our algorithm
discretizes the nonlinear Riccati equation satisfied by the 
derivative of the phase function via a Chebyshev spectral collocation method and 
applies the Newton-Kantorovich method to the resulting system of
nonlinear algebraic equations.
We prove that the iterates converge quadratically
to a nonoscillatory solution of the Riccati equation.
The quadratic convergence of the Newton-Kantorovich method and
the  simple form of the linearized equations
ensure that this procedure is extremely efficient.
Our algorithm  then extends the slowly varying phase function calculated
in the high-frequency regime throughout the solution domain 
 by solving a certain third order linear
ordinary differential equation related to the Riccati equation.
We describe the results of numerical experiments showing that
our algorithm is orders of magnitude faster than existing
schemes,  including 
the modified Magnus method \cite{ModMagnus},
the current state-of-the-art approach \cite{BremerPhase}
and the recently introduced  ARDC method \cite{ARDC}.

\end{abstract}

\author{Tara Stojimirovic}
\address{\vskip -1em Department of Mathematics, University of Toronto}

\author{James Bremer}
\ead{bremer@math.toronto.edu}
\address{\vskip -1em Department of Mathematics, University of Toronto\vskip -2.5em}

\end{frontmatter}

\begin{section}{Introduction}







Oscillatory second order linear ordinary differential equations are
widely used in scientific and numerical calculations.  
They arise, for instance, in plasma physics \cite{Hazeltine,Davidson,Agocs}, 
Hamiltonian dynamics \cite{2018_Pritula}, quantum mechanics \cite{Einaudi} and cosmology \cite{Jerome,Agocs}.  
Moreover,  numerous important families of  special functions,  such as the Jacobi polynomials, the 
Bessel functions and the spheroidal  wave functions, satisfy  such equations (see, for instance, \cite{DLMF}).

In many applications, the speed with which the differential equations
are solved is of great importance.  This is case for the 
calculations described in \cite{Agocs}, which entail the numerical solution
of billions of oscillatory ODEs, and there are numerous other such examples. 
Moreover, many calculations related to special
functions require a large number of function evaluations.
To give one example, it is necessary to rapidly evaluate
associated Legendre functions of many different orders and degrees
at a large number of points in order to apply the spherical harmonic transform 
via a butterfly algorithm \cite{BremerSHT}.

The principal result of this paper applies to equations of the form
\begin{equation}
y''(t) + \omega^2 q(t,\omega) y(t) = 0, \ \ \ -\infty < a < t < b < \infty,
\label{introduction:ode}
\end{equation}
where $q(t,\omega)$ is a smooth function  such that there exist positive constants $\qmin$, $\omega_0$  
and a positive integer $M$ with 
\begin{equation}
\left| \left(\frac{d}{dt}\right)^j q(t,\omega)\right| \leq 1
\ \ \ \mbox{for all}\ \ \ 0 \leq j \leq 2M
\label{introduction:qbound1}
\end{equation}
and
\begin{equation}
\qmin \leq  q(t , \omega)
\label{introduction:qbound2}
\end{equation}
for all $a \leq t \leq b$ and all $\omega \geq  \omega_0$.  
It is well known  that  under these assumptions,
there exists a basis $\{u,v\}$ in the space of solutions of (\ref{introduction:ode}) such that 
\begin{equation}
u(t) =  
\frac{\exp\left(i \omega \int_a^t \sqrt{q(s,\omega)}\, ds\right)}{\sqrt{\omega}\, q(t,\omega)^{\frac{1}{4}}} 
+ \mathcal{O}\left(\frac{1}{\omega}\right)
\ \ \ \mbox{as} \ \ \ \omega\to\infty
\label{introduction:lg}
\end{equation}
and
\begin{equation}
v(t) =  
\frac{\exp\left(-i \omega \int_a^t \sqrt{q(s,\omega)}\, ds\right)}{\sqrt{\omega}\, q(t,\omega)^{\frac{1}{4}}} 
+ \mathcal{O}\left(\frac{1}{\omega}\right)
\ \ \ \mbox{as} \ \ \ \omega\to\infty.
\label{introduction:lg2}
\end{equation}
These expressions are known as Liouville-Green or first-order WKB approximates for (\ref{introduction:ode}),
and they are usually derived by analyzing the Riccati equation
\begin{equation}
r'(t) + (r(t))^2 + \omega^2 q(t,\omega) = 0
\label{introduction:riccati}
\end{equation}
satisfied by the logarithmic derivatives of solutions of (\ref{introduction:ode}).
Throughout this paper, we will refer to any function $\psi$ such that $\exp(\psi(t))$ is a solution
of (\ref{introduction:ode}) as an {\it exponential phase function}
for (\ref{introduction:ode}), so that  the solutions of the Riccati equation
(\ref{introduction:riccati}) are the derivatives of exponential phase functions.
Among other things, (\ref{introduction:lg}) and (\ref{introduction:lg2}) imply that
the solutions of (\ref{introduction:ode}) are oscillatory, and that the parameter $\omega$
is a reasonable proxy for their frequency of their oscillations.
It follows that when a standard ODE solver is applied to (\ref{introduction:ode}), its running
time grows roughly linearly with $\omega$.

Here, we give an elementary proof that under the assumptions (\ref{introduction:qbound1})
and (\ref{introduction:qbound2}), there exists a nonoscillatory solution of the Riccati equation
satisfied by the logarithmic derivatives of the solutions of (\ref{introduction:ode}).   The solution
is nonoscillatory in the sense that it can be approximated to a fixed accuracy via
a Chebyshev expansion the number of terms of which is independent of the frequency
parameter  $\omega$, at least for sufficiently large $\omega$.
Moreover, we show that when the Riccati equation is discretized using a standard Chebyshev spectral collocation method
and the Newton-Kantorovich method is applied to the resulting system of nonlinear algebraic
equations, the iterates converge quadratically to a vector giving the values
of this nonoscillatory phase function at the collocation nodes,
again provided that the frequency parameter $\omega$ is sufficiently large.
Because of the quadratic convergence of the iterates and 
the fact that the linearized equations which arise 
are of an extremely simple type that can be inverted rapidly,  this procedure is considerable faster
than existing methods for constructing nonoscillatory phase functions
representing solutions of ordinary differential equations of the form (\ref{introduction:ode}).

We go on to describe a numerical algorithm for solving a large class of oscillatory differential equations.
In particular, our algorithm applies in cases when the solutions are highly oscillatory in some regions and slowly varying in others.
It runs in time independent of frequency,
and obtains accuracy on the order of the condition number of evaluation of the solutions.
Our method operates by constructing a phase function which is represented via a piecewise Chebyshev expansion.
More explicitly, it adaptively subdivides the solution domain $[a,b]$ of the equation into a collection
of discretization subintervals $[a_1,b_1],\ldots,[a_m,b_m]$ and, on each subinterval, the phase
function is represented via a Chebyshev expansion of a fixed
order $k$.
The Newton-Kantorovich method is used to construct the polynomial expansion over each ``high-frequency interval.''
The linearized equations which arise during the course of the Newton-Kantorovich iterations
are not uniquely solvable in the low-frequency regime, however, and
so our algorithm extends the  phase function into the rest of the solution domain by 
solving initial and terminal value problems for  Appell's equation.  Appell's equation is a certain third order linear
ordinary differential equation related to the Riccati equation (\ref{introduction:riccati}).
We prefer Appell's equation to the Riccati equation for two reasons.  First, it is linear and this means
that it can be solved more rapidly than the Riccati equation, at least in the low-frequency regime
where the simple form of the linearized Riccati equation is difficult to  exploit.
Secondly, Riccati's equation becomes numerically degenerate when the coefficient $q$ is zero or close to it,
whereas,  as observed in \cite{BremerPhase2}, Appell's equation can be solved in a numerical
stable fashion even when $q$ is close to $0$ or  negative and of large magnitude. 
As described here, our algorithm is limited to equations whose coefficients are nonnegative
throughout the solution domain; however, with minor alterations, it could be applied in cases in which 
the coefficient becomes negative on some parts of its domain.  
Over such regions, the solutions of the differential equation behave like combinations of rapidly increasing and 
decreasing exponential functions, but they can nonetheless be accurately represented via phase functions (see \cite{BremerPhase2}).

Almost all specialized techniques for solving oscillatory ordinary differential equations
are based  on the observation that such equations admit nonoscillatory exponential representations.
One of the most widely used algorithms is the modified Magnus expansion method described in \cite{ModMagnus}.
It exploits the existence of efficient exponential representations of the solutions
of oscillatory differential equations by using an exponential integrator as the basis of a step method.
Moreover, in each step, the equation is preconditioned by the solution of the constant
coefficient equation obtained by freezing the coefficients at the midpoint.
The running time of the modified Magnus method grows 
as $\mathcal{O}\left(\omega^{3/4}\right)$ when it is applied to a scalar equation of the form
(\ref{introduction:ode}).  It should be noted, though, that it applies to a 
wider class of differential equations than that considered here,
including quite general systems of linear ordinary differential equations.
Like the modified Magnus method, the scheme of  \cite{Lubich} applies to a large class of systems of 
differential equations.  It operates by  constructing a preconditioner using the eigendecomposition 
of the system's coefficient matrix.  When applied to scalar equations of the type appearing here, 
its running time  grows as $\mathcal{O}\left(\sqrt{\omega}\right)$.
The WKB marching method of \cite{Arnold,Korner},
which is specific to second order linear ordinary differential equations,
uses  second order WKB expansions as a preconditioner in regions 
where the solutions are rapidly oscillating and applies  a standard Runge-Kutta method in the 
low-frequency regime.    The running time of the WKB marching method also grows as  
$\mathcal{O}\left(\sqrt{\omega}\right)$ in typical cases.

Another class of numerical methods uses asymptotic approximations to represent
solutions of  second order linear ordinary differential equations
in the high-frequency regime directly, as opposed to using asymptotic
approximations as preconditioners. The article \cite{Agocs} introduces a solver for 
second order linear  ordinary differential equations that
 represents  solutions via low-order WKB expansions  in regions where
the solutions are  highly oscillatory,  and applies a standard Runge-Kutta method
in the low-frequency regime.           The running time of this algorithm is independent 
of frequency in certain special cases, for instance
when the solutions oscillate at extremely high frequency over the whole solution domain,
but, in the general case, it grows as $\mathcal{O}\left(\omega\right)$.
The ARDC scheme \cite{ARDC} improves on \cite{Agocs} by using higher order WKB-like 
approximations constructed numerically through an iterative method
to represent solutions in the oscillatory regime. 
The running time of \cite{ARDC} is independent of frequency, but
it is significantly slower and no more widely applicable than the earlier frequency-independent method introduced in 
\cite{BremerPhase} (see Section~\ref{section:experiments2}).

The iterative method used by  \cite{ARDC} is strongly related to the approach of 
\cite{SpiglerPhase2,SpiglerZeros,SpiglerPhase1}.    There,
symbolic formulas  which represent asymptotic approximations
of solutions of the differential equation satisfied by the square of the imaginary 
part of solutions of the Riccati equation are constructed via an iterative scheme.
The calculations are conducted symbolically rather than numerically because  each iteration requires a further
derivative of the coefficient.  
In \cite{SheehanOlver1}, a method for solving a large class of oscillatory differential
equations in time independent of frequency is described.  It applies the 
 ``differential GMRES'' procedure described in \cite{SheehanOlver2} to a auxiliary
system derived from the differential equation under consideration.
It also uses symbolic calculations because it requires high order derivatives of the entries
of the coefficient matrix.

The method of \cite{BremerPhase} appears to be the current
state-of-the-art algorithm for numerical solution of the class of 
equations considered here.
It operates by solving the equation
\begin{equation}
\left(\alpha'(t)\right)^2 - \omega^2 q(\omega,t) 
- \frac{3}{4}\left(\frac{\alpha''(t)}{\alpha'(t)}\right)^2 
+ \frac{1}{2} \frac{\alpha'''(t)}{\alpha'(t)} = 0
\label{introduction:kummer}
\end{equation}
obtained by considering the real and imaginary parts of (\ref{introduction:riccati})
separately in order to calculate a function $\alpha$ such that
\begin{equation}
\frac{\cos\left(\alpha(t)\right)}{\sqrt{\alpha'(t)}}
\ \ \ \mbox{and}\ \ \
\frac{\sin\left(\alpha(t)\right)}{\sqrt{\alpha'(t)}}
\label{introduction:uv}
\end{equation}
constitute a basis in the space of solutions of (\ref{introduction:ode}).
We refer any $\alpha$ such that the function $u$ and $v$ defined
in (\ref{introduction:uv}) are a pair of independent solutions of (\ref{introduction:ode})
as a {\it trigonometric phase function} for (\ref{introduction:ode}).
Equation~(\ref{introduction:kummer}) is satisfied
by the derivatives of trigonometric phase functions,
and we refer to it as Kummer's equation
after E.~E.~Kummer who studied it in \cite{Kummer}.
The algorithm of \cite{BremerPhase} runs in time 
independent of frequency and like the method of this paper, but unlike those proposed in 
\cite{ARDC,Arnold,Korner,Agocs}, the same exponential representations
of solutions are used throughout the solution domain.   This is highly conducive to certain 
calculations, such as the computation of the  zeros of solutions \cite{BremerZeros} and the
application of Strum-Liouville eigentransforms
\cite{BremerSHT}.
Almost all of the solutions of the Riccati equation are oscillatory or singular
and the main focus of \cite{BremerPhase} is  a   mechanism which selects the desired 
nonoscillatory  solution.     This is done by smoothly deforming $q$
so that it  is  constant on $[a,(3a+b)/4]$ and equal  to the original coefficient $q$ on  $[(a+3b)/4,b]$.
The value of an appropriate nonoscillatory solution of the Riccati equation
corresponding to the deformed coefficient  is known
at $a$, and by solving an initial value problem for the modified Riccati equation,
the value at $b$ of a slowly-varying solution of the original Riccati
equation corresponding to (\ref{introduction:ode}) is obtained.  
The derivative of the desired phase function is then constructed by solving a 
terminal value problem for the original Riccati equation.  One of the main observations
of this paper is that the machinery of \cite{BremerPhase} is unnecessary
 because   discretizing the Riccati equation in one of the simplest
ways possible and applying one of the most basic strategies for solving 
the resulting nonlinear system yields a faster method.

The remainder of this paper is structured as follows.  Section~\ref{section:preliminaries}
reviews the necessary  mathematical and numerical preliminaries.   Section~\ref{section:nonoscillatory}
gives an elementary proof of the existence of a nonoscillatory solution of the Riccati equation.
In Section~\ref{section:riccati}, we prove that when the Riccati equation is discretized
via a Chebyshev spectral method and the Newton-Kantorovich method method is applied to the
resulting system of nonlinear algebraic equations, the iterates converge to a vector which represents
the nonoscillatory solution whose existence is shown in Section~\ref{section:nonoscillatory}.
Section~\ref{section:algorithm} details
our numerical algorithm, and we go on to describe the results of numerical experiments
conducted to demonstrate its properties in Section~\ref{section:experiments}. We close with 
a few brief remarks regarding this work and future directions for research in Section~\ref{section:conclusion}.

\end{section}

\begin{section}{Mathematical and Numerical Preliminaries}
\label{section:preliminaries}

%
%

\begin{subsection}{Notation and conventions}

We denote the Fr\'echet derivative of a map $F:X \to Y$ between Banach
spaces at the point $x$ via $F'(x)$.  Assuming it exists, 
the Fr\'echet derivative is an element of the space of bounded linear
functions $X \to Y$, for which we use the notation  $L(X,Y)$.
We will use $B_r(x)$ for the ball of radius $r>0$ centered at the point $x \in  X$.
A function $F:\Omega \subset X \to Y$ given on an open subset $\Omega$ of a Banach space
$X$ is said to be continuously differentiable on $\Omega$ 
provided the map $\Omega \to L(X,Y)$ given by $x \mapsto F'(x)$ 
is continuous.

In Section~\ref{section:riccati}, we will use
vectors in $\mathbb{R}^n$ to represent certain functions and it is
convenient to introduce a convention for distinguishing between a function and the vector
representing it.  Accordingly, throughout this paper we will display the names of vectors 
in the space $\mathbb{R}^n$
using a bold font.   
Moreover, we let  $\mbox{diag}\left(\mathbf{x}\right)$ denote the $n\times  n$ diagonal matrix 
whose diagonal entries are the elements of the vector $\mathbf{x} \in \mathbb{R}^n$.
We use the notation $\mathbf{v} \circ \mathbf{w}$ 
to denote the entrywise or Hadamard product of the vectors $\mathbf{v}$ and $\mathbf{w}$ --- that is,
the vector whose $j^{th}$ entry is the product of the $j^{th}$ entry of the vector
$\mathbf{v}$ with the $j^{th}$ entry of  the vector $\mathbf{w}$.


We use $T_j$ for the Chebyshev polynomial of degree $j$.
We denote nodes of the  $k$-point Chebyshev extremal grid on $[-1,1]$ by 
\begin{equation}
\tcheb{1}, \tcheb{2}, \ldots, \tcheb{k},
\end{equation}
so that 
\begin{equation}
\tcheb{i} =  \cos\left(\pi \frac{k-i}{k-1}\right).
\end{equation}
We use $\mathscr{D}_k$ for the $k\times k$ spectral differentiation matrix
which takes the vector
\begin{equation}
\left(
\begin{array}{cccccc}
p\left(\tcheb{1}\right) & p\left(\tcheb{2}\right) & \cdots & p\left(\tcheb{k}\right)
\end{array}
\right)
\label{preliminaries:pvals}
\end{equation}
of values of a polynomial $p$ of degree less than $k$ 
at the Chebyshev nodes to the vector
\begin{equation}
\left(
\begin{array}{cccccc}
p'\left(\tcheb{1}\right) & p'\left(\tcheb{2}\right) & \cdots & p'\left(\tcheb{k}\right)
\end{array}
\right)
\end{equation}
of the values of its derivatives at the same nodes.
Moreover, we let  $\mathscr{I}_k$ be the spectral integration matrix which takes the vector
of values (\ref{preliminaries:pvals}) to the vector  
\begin{equation}
\left(
\begin{array}{cccccc}
q\left(\tcheb{1}\right) & q\left(\tcheb{2}\right) & \cdots &
q\left(\tcheb{k}\right)
\end{array}
\right),
\label{preliminaries:qvals}
\end{equation}
where 
\begin{equation}
q(t) = \int_{-1}^t p(s)\ ds.
\end{equation}
Finally, we let $\mathscr{C}_k$ denote the matrix which takes the vector
of values (\ref{preliminaries:pvals}) of a polynomial $p$ of degree less than
$k$
to the vector 
\begin{equation}
\left(
\begin{array}{cccccc}
a_0 & a_1 & \cdots & a_{k-1}
\end{array}
\right),
\end{equation}
of coefficients in the Chebyshev expansion
\begin{equation}
p(t) = \sum_{j=0}^{k-1} a_j T_j(t).
\end{equation}

\end{subsection}

%
%

\begin{subsection}{The Newton-Kantorovich theorem}

In \cite{Kantorovich},  Newton's method is generalized to the case of maps between Banach spaces and
conditions for its convergence are given.  Here, we state a simplified version of the 
Newton-Kantorovich theorem,  which can be found as Theorem~7.7-4 in
Section~7.7 of \cite{Ciarlet}.   

\begin{theorem}[Newton-Kantorovich]
\label{theorem:nk}
Suppose that $\Omega$ is an open subset of the Banach space $X$, 
that $Y$ is a Banach space and that $F: \Omega \subset X \to Y$ is continuously differentiable.
Suppose also that there exist a point $x_0 \in \Omega$ and
constants $\lambda$ and $\eta$  such that
\begin{enumerate}

\item
$\begin{aligned}
F'(x_0)
\end{aligned}$
admits an inverse $F'(x_0)^{-1} \in L(Y,X)$,
\vskip .25em

\item 
$\begin{aligned}
B_\eta\left(x_0\right) \subset \Omega,
\end{aligned}$
\vskip .25em

\item
$\begin{aligned}
0 < \lambda < \frac{\eta}{2},
\end{aligned}$
\vskip .25em

\item
$\begin{aligned}
\left\| F'(x_0)^{-1} F(x_0) \right\| \leq \lambda \ \ \ \mbox{and}
\end{aligned}$
\vskip .25em

\item
$\begin{aligned}
\left\|F'(x_0)^{-1} \left(F'(x) - F'(y) \right) \right\| \leq \frac{1}{\eta} \left\|x-y\right\|
\  \mbox{for all}  \  x,y \in B_\eta(x_0).
\end{aligned}$
\end{enumerate}
Then, $F'(x)$ has a bounded inverse $F'(x)^{-1} \in L(Y,X)$ for each $x \in B_\eta(x_0)$, 
the sequence $\{x_n\}$ defined by
\begin{equation}
x_{n+1} = x_n - F'(x_n)^{-1} F(x_n)
\end{equation}
is contained in the open ball
\begin{equation}
 B_{\eta^-}(x_0), \ \ \ \mbox{where} \ \ \ \eta^- = \eta \left(1 - \sqrt{1-\frac{2\lambda}{\eta}}\right) \leq \eta,
\end{equation}
and $\{x_n\}$ converges to a zero $x^*$ of $F$.
Moreover,  $x^*$ is the only  zero of $F$ in the ball $B_\eta(x_0)$
and, for each $n \geq 0$, we have
\begin{equation}
\left\| x_n - x^* \right\| \leq \frac{\eta}{2^n} \left( \frac{\eta^-}{\eta}\right)^{2^n}.
\end{equation}

\end{theorem}

We will refer to the equation
\begin{equation}
F'\left(x_n\right)h = - F\left(x_n\right),
\end{equation}
which arises when calculating the $(n+1)^{st}$ Newton-Kantorovich iterate as the linearization of the 
equation $F\left(x\right) = 0$ around the point $x_n$.

\end{subsection}

%
%

\begin{subsection}{The Banach fixed point theorem}
\label{section:fixed}

A proof of the following well-known theorem can be found in Section~3.7 of \cite{Ciarlet},
among many other sources.

\begin{theorem}[Banach fixed point theorem]
\label{theorem:fp}
Suppose that $F:X \to X$ is a map between Banach spaces, and that
there exists a constant $0 \leq \gamma < 1$ such that
$\left\|F(x_1)-F(x_2)\right\| \leq \gamma \left\|x_1-x_2\right\|$
for all $x_1,x_2 \in X$.  Then $F$ has a unique
fixed point $x^*$.  Moreover, given any $x_0 \in X$, the sequence
$\{x_n\}$ defined via $x_{n+1} = F( x_n)$ converges to $x^*$, and
the following estimate holds
\begin{equation}
\left\|x_n-x^*\right\| \leq \frac{\gamma^n \left\| F(x_0)-x_0\right\| }{1-\gamma}.
\end{equation}
\end{theorem}

An immediate consequence of Theorem~\ref{theorem:fp} is that a system of linear algebraic equations
of the form
\begin{equation}
(I + S ) \mathbf{x} = \mathbf{y}
\label{preliminaries:pid}
\end{equation}
can be solved extremely rapidly via the iteration
\begin{equation}
\mathbf{x_0} = 0 \ \ \ \mbox{and} \ \ \ \mathbf{x_{n+1}} = -S \mathbf{x_n} + \mathbf{y}
\label{preliminaries:fixed}
\end{equation}
as long as there is some operator norm $\| \cdot \|_0$ such that $\|S\|_0 \ll 1$.
The spectral radius $\rho(A)$ of a matrix $A$ is the maximum of the magnitudes of its eigenvalues.
It is well known (see, for instance, \cite{Horn}) 
that given any $\epsilon > 0$, there exists a matrix norm $\|\cdot\|_\epsilon$ such that 
\begin{equation}
\rho(S) \leq \|S\|_\epsilon \leq \rho(S)+\epsilon,
\end{equation}
so the iteration (\ref{preliminaries:fixed}) converges rapidly to a solution
of (\ref{preliminaries:pid}) assuming $\rho(S) \ll 1$.

\end{subsection}

%
%

\begin{subsection}{The Riccati equation and its variants}

In this section, we briefly discuss the relationship between the solutions
of the ordinary differential equations (\ref{introduction:ode}), 
(\ref{introduction:riccati}) and (\ref{introduction:kummer}),
and introduce a related third order linear ordinary differential equation.

It can be verified through direct computation that if
\begin{equation}
y(t) = \exp\left(\int r(s)\, ds\right)
\end{equation}
satisfies (\ref{introduction:ode}), then $r$ is a solution of the Riccati equation (\ref{introduction:riccati}).
It follows by separating out the real and imaginary parts of the Riccati equation
that its solutions are of the form
\begin{equation}
r(t) = i \alpha'(t) - \frac{1}{2}\frac{\alpha''(t)}{\alpha'(t)},
\label{variants:ralpha}
\end{equation}
where $\alpha$ satisfies Kummer's equation (\ref{introduction:kummer}).  In particular,
the functions
\begin{equation}
u(t) = \frac{\sin\left(\alpha(t)\right)}{\sqrt{\alpha'(t)}}
\ \ \  \mbox{and}\ \ \ 
v(t) = \frac{\cos\left(\alpha(t)\right)}{\sqrt{\alpha'(t)}}
\label{variants:uv}
\end{equation}
constitute a basis of real-valued functions in the space of solutions of (\ref{introduction:ode}).
It is clear from (\ref{variants:uv}) that the sum of squares of $u$ and $v$ is the  reciprocal
of $\alpha'$, and it can be shown that if  $u$, $v$ is any pair of real-valued solutions of (\ref{introduction:ode}) 
whose Wronskian is $1$, then 
\begin{equation}
\alpha'(t) = \frac{1}{(u(t))^2 + (v(t))^2}
\end{equation}
is a solution of Kummer's equation. 
A somewhat lengthy, but straightforward, computation
shows that the modulus function
\begin{equation}
m(t) = (u(t))^2 + (v(t))^2 = \frac{1}{\alpha'(t)}
\end{equation}
is a solution of the third order linear ordinary differential equation
\begin{equation}
m'''(t) + 4 \omega^2 q(\omega,t) m'(t) + 2 \omega^2 \frac{dq}{dt} (\omega,t) m(t) = 0,
\label{variants:appell}
\end{equation}
which we refer to as Appell's equation in light of \cite{appell}.

In this way, each ``phase function'' for (\ref{introduction:ode})
corresponds with a solution of (\ref{introduction:riccati}), a solution of (\ref{introduction:kummer}),
a modulus function satisfying (\ref{variants:appell}) and  
a pair of real-valued solutions of (\ref{introduction:ode}) whose Wronskian is $1$.
For the sake of clarity and precision, we refer to functions $\psi$ such that
$\exp(\psi(t))$ is a solution of (\ref{introduction:ode})
as  exponential phase functions, and we
say that $\alpha$ is a trigonometric phase function if the functions
defined in  (\ref{variants:uv}) constitute
a basis in the space of solutions of (\ref{introduction:ode}).
Given the close relationship (\ref{variants:ralpha}) between the solutions of the Riccati equation
and those of Kummer's equation, however, there is little real distinction
between exponential and trigonometric phase functions.
 
\end{subsection}

%
%

\begin{subsection}{A bound on the magnitude of Chebyshev coefficients}

The Chebyshev expansion of a function $f:[a,b] \to \mathbb{C}$ is the series
\begin{equation}
\sum_{n=0}^\infty a_n T_n\left(\frac{2}{b-a} t - \frac{b+a}{b-a} \right),
\label{chebyshev:series}
\end{equation}
where
\begin{equation}
a_0 = \frac{1}{\pi} \int_{-1}^1 f\left(\frac{b-a}{2} t + \frac{b+a}{2} \right)\,  \frac{dt}{\sqrt{1-t^2} }
\label{chebyshev:coefs0}
\end{equation}
and
\begin{equation}
a_n = \frac{2}{\pi} \int_{-1}^1 f\left(\frac{b-a}{2} t + \frac{b+a}{2} \right) T_n\left(t\right)\,  \frac{dt}{\sqrt{1-t^2} }
\ \ \ \mbox{for all}\ \ \ n >0.
\label{chebyshev:coefs}
\end{equation}
%
The following theorem, which can be found in a slightly different form in \cite{trefethen}, gives an upper bound
on the magnitudes of the Chebyshev coefficients of $f$.
\begin{theorem}
\label{theorem:cheb}
If $f,f',f'',\ldots,f^{(k-1)}$ are absolutely continuous on $[a,b]$ and 
\begin{equation}
V =  \left(\frac{b-a}{2}\right)^k \int_{-1}^1 \left|\frac{f^{(k)}\left(\frac{b-a}{2} t + \frac{b+a}{2} \right)}{\sqrt{1-t^2}}\right|\, dt,
\end{equation}
then we have the following bound on the magnitudes of the Chebyshev coefficients $a_n$ of $f$:
\begin{equation}
\left|a_n\right| \leq \frac{2V}{\pi n (n-1)(n-2)\cdots(n-k)}  \ \ \ \mbox{for all} \ \ \ n > k.
\end{equation}
\end{theorem}

\end{subsection}

\end{section}

 \begin{section}{Existence of Nonoscillatory Phase Functions}
\label{section:nonoscillatory}

In this section, we first give an elementary argument showing that under the assumptions
(\ref{introduction:qbound1}) and (\ref{introduction:qbound2}), the Riccati equation 
admits a solution which is  nonoscillatory in the sense that it can be approximated via a  
polynomial expansion whose number of terms $k$  is independent
of $\omega$, at least for large values of $\omega$.
In many cases of interest, $q$ is nonoscillatory in more robust sense.
Accordingly, we go on to discuss the principal 
theorem of \cite{BremerPhase}, which 
shows the existence of a  phase function which is nonoscillatory in a stronger sense
under regularity conditions on $q$ which are more stringent
than (\ref{introduction:qbound1}) and (\ref{introduction:qbound2}).

\begin{subsection}{Elementary proof of the existence of a nonoscillatory phase function }
\label{section:nonoscillatory1}

We let  $r_M$ be defined via the formula
\begin{equation}
r_M(t) = \omega u_0(t) + u_1(t) + \cdots + u_{M}(t) \omega^{1-M},
\label{nonoscillatory:rm}
\end{equation}
where
\begin{equation}
u_0(t) = i \omega \sqrt{q(t,\omega)}, \ \ \ u_1(t) = -\frac{1}{4}\frac{q'(t)}{q(t)}
\end{equation}
and
\begin{equation}
u_{n}(t) = -\frac{1}{2 u_0(t)} \left(u_{n-1}'(t) + \sum_{j=1}^{n-1} u_j(t) u_{n-j}(t) \right)
\ \ \ \ \mbox{for all} \ \ \ j \geq 1.
\end{equation}
We will first show that there exists a solution $r$ of (\ref{introduction:riccati})
such that 
\begin{equation}
\left\|r - r_M\right\|_{L^\infty\left(\left[a,b\right]\right)} =  \mathcal{O}\left(\frac{1}{\omega^M}\right) \ \ \ \mbox{as}\ \ \ \omega\to\infty.
\label{nonoscillatory:rbound}
\end{equation}
%
%
To that end, we insert the expression
\begin{equation}
r(t) = r_M(t) + \delta(t)
\label{nonoscillatory:delta}
\end{equation}
into (\ref{introduction:riccati}), which yields the equation
\begin{equation}
\delta'(t) + 2 r_M(t)  \delta(t) + \left(\delta(t)\right)^2 = - G(t),
\label{nonoscillatory:delta_eq}
\end{equation}
where
\begin{equation}
G(t)=r_M'(t) + \left(r_M(t)\right)^2 + \omega^2 q(t,\omega).
\end{equation}
We observe that  the assumptions (\ref{introduction:qbound1}) and (\ref{introduction:qbound2})
together with the fact that $u_j(t)$ involves only $q$ and the first $j$ derivatives of $q$,
imply 
\begin{equation}
\left\|G\right\|_{L^\infty\left(\left[a,b\right]\right)} = \mathcal{O}\left(\frac{1}{\omega^{M-1}}\right) 
\ \ \ \mbox{and}\ \ \ 
\left\|G'\right\|_{L^\infty\left(\left[a,b\right]\right)} = \mathcal{O}\left(\frac{1}{\omega^{M-1}}\right)  \ \ \ \mbox{as}\ \ \ \omega\to\infty.
\label{nonoscillatory:Fbounds}
\end{equation}
We proceed by setting
\begin{equation}
H(t) = 2\int_a^t r_M(s),
\end{equation}
and multiplying both sides of (\ref{nonoscillatory:delta_eq}) by $\exp(H(t))$ to obtain 
\begin{equation}
\exp(H(t))\delta'(t) + 2 \exp(H(t)) r_M(t) \delta(t) =- \exp(H(t)) \left(\delta(t)\right)^2  +  \exp(H(t))G(t).
\label{nonoscillatory:delta_eq2}
\end{equation}
Since the left-hand side of (\ref{nonoscillatory:delta_eq2}) is the derivative of $\delta(t) \exp(H(t))$,
we can reformulate it as the integral equation
\begin{equation}
 \delta(t)  = -\exp(-H(t)) \int_{a}^t \exp(H(s)) \left(  \left(\delta(s)\right)^2  + G(s)\right)\, ds.
\label{nonoscillatory:inteq}
\end{equation}
We will show that when $\omega$ is sufficiently large, the operator
\begin{equation}
S\left[\delta\right](t) = -\exp(-H(t)) \int_{a}^t \exp(H(s)) \left(  \left(\delta(s)\right)^2  + G(s)\right)\, ds
\end{equation}
is a contraction on a closed ball $B$ of whose radius $2C/\omega^M$ centered at $0$ in  $L^\infty\left([a,b]\right)$,
where $C$ is an appropriately chosen constant.  It will then follow from Theorem~\ref{theorem:fp} that there is 
a solution of  (\ref{nonoscillatory:inteq}) in this ball, and the desired bound (\ref{nonoscillatory:rbound})
is a consequence of the existence of this solution and (\ref{nonoscillatory:delta}).

Our assumptions on $q$ and its derivatives
imply that the functions $\left|\exp(\pm H(t))\right|$ are bounded independent of $\omega$.
It follows that there exists a constant $C>0$ which is independent of $\omega$ and such that
\begin{equation}
\left\| -\exp(-H(t)) \int_{a}^t \exp(H(s)) f(s)\, ds \right\|_{L^\infty\left([a,b]\right)}
\leq C \| f \|_{L^\infty\left([a,b]\right)}
\label{nonoscillatory:bound1}
\end{equation}
for any $f \in L^\infty\left(\left[a,b\right]\right)$.  In particular, we have
\begin{equation}
\left\| -\exp(-H(t)) \int_{a}^t \exp(H(s)) \left(\delta(s)\right)^2\, ds \right\|_{L^\infty\left([a,b]\right)}
\leq C \| \delta \|_{L^\infty\left([a,b]\right)}^2.
\label{nonoscillatory:bound1.5}
\end{equation}
Now, integration by parts shows that 
\begin{equation}
\begin{aligned}
 &\exp(-H(t)) \int_{a}^t \exp(H(s)) G(s)\, ds \\
= 
&\exp(-H(t)) \int_{a}^t \exp(H(s))H'(s) \frac{G(s)}{H'(s)}\, ds \\
=
 &\exp(-H(t))\left( \exp(H(t))\frac{G(t)}{H'(t)} -\exp(H(a))\frac{G(a)}{H'(a)}
\right. \\ \left.
\right.&\left.-\int_a^t \exp(H(s)) \frac{G'(s)}{H'(s)}\, ds
+\int_a^t \exp(H(s)) \frac{G(s)H''(s)}{(H'(s))^2}\, ds
 \right).
\end{aligned}
\label{nonoscillatory:intparts}
\end{equation}
It is easy to see that 
\begin{equation}
\frac{1}{H'(s)} = \mathcal{O}\left(\frac{1}{\omega}\right)
\ \ \ \mbox{and} \ \  \
\frac{H''(s)}{(H'(s))^2} =  \mathcal{O}\left(\frac{1}{\omega}\right),
\label{nonoscillatory:Hbounds}
\end{equation}
and it follows  from (\ref{nonoscillatory:Hbounds}), the fact that the functions $\exp(\pm H(t))$ are bounded
independent of $\omega$, (\ref{nonoscillatory:intparts})
and  (\ref{nonoscillatory:Fbounds}) that  we can adjust the constant $C$ so that
\begin{equation}
\begin{aligned}
\left\| \exp(-H(t)) \int_{a}^t \exp(H(s)) G(s)\, ds\right\|_{L^\infty\left([a,b]\right)}
\leq \frac{C}{\omega^M}
\end{aligned}
\label{nonoscillatory:bound2}
\end{equation}
holds for all sufficiently large $\omega$.
Together (\ref{nonoscillatory:bound1.5}) and (\ref{nonoscillatory:bound2}) imply that
\begin{equation}
\left\|S\left[\delta\right]\right\|_{L^\infty\left([a,b]\right)}
\leq  
\frac{4 C^3}{\omega^{2M}} + \frac{C}{\omega^M}
\end{equation}
whenever $\|\delta\|_{L^\infty\left([a,b]\right)} \leq 2 C /\omega^M$.
It follows that for all sufficiently large $\omega$ and all $\delta$
in the ball $B$ of radius $2C/\omega^M$ centered at $0$,
\begin{equation}
\left\|S\left[\delta\right]\right\|_{L^\infty\left([a,b]\right)} \leq \frac{2C}{\omega^M};
\end{equation}
that is, the operator $S$ preserves the ball $B$.  Next, we observe that for all $\delta$ in $B$, we have
\begin{equation}
\begin{aligned}
&\left|\exp(-H(t)) \int_a^t \exp(H(s)) \left( \delta_1(s) - \delta_2(s) \right)\left( \delta_1(s) + \delta_2(s) \right)\  ds\right|\\
\leq
&C \left\| \delta_1 - \delta_2 \right\|_{L^\infty\left([a,b]\right)} \left\| \delta_1 + \delta_2 \right\|_{L^\infty\left([a,b]\right)}
\\
\leq &\frac{4C^2}{\omega^M} \| \delta_1-\delta_2\|_{L^\infty\left([a,b]\right)},
\end{aligned}
\end{equation}
from which we see that $S$ is a contraction on the ball $B$ provided $\omega$ is sufficiently large.
By applying Theorem~\ref{theorem:fp}, we conclude that $S$ admits a unique fixed point in $B$ and, as previously
discussed, the desired bound (\ref{nonoscillatory:rbound}) follows.

We now apply Theorem~\ref{theorem:cheb} to the function $r_M$.  
Since the first $M$ derivatives of $r_M$ depend on the first $2M$ derivatives of $q$,
(\ref{introduction:qbound1}) and (\ref{introduction:qbound2}) imply that 
\begin{equation}
V =  \left(\frac{b-a}{2}\right)^M \int_{a}^b \left|\frac{r_M^{(M)}\left(\frac{b-a}{2} t + \frac{b+a}{2} \right)}{\sqrt{1-t^2}}\right|\, dt,
\end{equation}
where $r_M^{(M)}$ is the $M^{th}$ derivative of the function $r_M$,
is bounded independent of $\omega$.  We invoke Theorem~\ref{theorem:cheb} to see that
\begin{equation}
r_M(t) = \sum_{n=0}^\infty a_n T_n\left(\frac{b-a}{2} t + \frac{b+a}{2}\right),
\end{equation}
where
\begin{equation}
\left|a_n\right| \leq \frac{V}{\pi n(n-1)\cdots(n-M)} \ \ \ \mbox{for all} \ \ n > M.
\end{equation}
It follows that there exists a constant $D$ independent of $\omega$ such that
\begin{equation}
\left\| r_M(t)  - \sum_{0}^{k-1} a_n T_n\left(\frac{b-a}{2} t + \frac{b+a}{2}\right) \right\|_{L^\infty\left(\left[a,b\right]\right)}
\leq \frac{D}{k^M}.
\label{nonoscillatory:polybound}
\end{equation}
Combining (\ref{nonoscillatory:polybound}) and (\ref{nonoscillatory:rbound}) shows that given $\epsilon >0$ we can choose an integer
$k$ and a solution $r$ of (\ref{introduction:riccati}) such that
\begin{equation}
\left\| r(t)  - \sum_{0}^{k-1} a_n T_n\left(\frac{b-a}{2} t + \frac{b+a}{2}\right) \right\|_{L^\infty\left(\left[a,b\right]\right)}
< \epsilon
\end{equation}
for all sufficiently large $\omega$. Moreover, the proceeding analysis implies
 that the necessary values of $k$ and $\omega$ are on the
order of  $\epsilon^{-\frac{1}{M}}$.

\end{subsection}

\begin{subsection}{A stronger results regarding nonoscillatory solutions of the Riccati equation}

We now reproduce a theorem proved in \cite{BremerRokhlin}
which gives conditions under which the Riccati equation admits solutions
that are nonoscillatory in a strong sense.    
It applies to equations given on the 
entire real line.  However, it is easy to transform to an equation of the type
(\ref{introduction:ode}) to one given on the real line;  for instance, by letting
\begin{equation}
t = \frac{b-a}{2} \tanh\left(x\right)+ \frac{b+a}{2}.
\end{equation}
Moreover, in order to simplify notation, we will restrict attention to
the case in which the coefficient $q$ does not depend on $\omega$.  The theorem 
can be  applied in the general case, assuming  that $q$ satisfies the necessary  conditions
uniformly in $\omega$, for sufficiently large $\omega$.
We let 
\begin{equation}
p(t) = \frac{1}{q(t)}  \left(\frac{5}{4} \left(\frac{q'(t)}{q(t)}\right)^2 - \frac{q''(t)}{q(t)}\right)
\end{equation}
and define the new variable $x$ via 
\begin{equation}
x(t) =  \int_0^t \sqrt{q(s)}\, ds.
\label{preliminaries:xt}
\end{equation}
The theorem requires a condition on the function $p(t(x))$.  At first glance, it might
seem that the relationship between this function and the coefficient $q$ is quite complicated;
however, when $p$ is viewed as a function of $x$, it is simply equal to a constant multiple of  
the Schwarzian derivative $\{t,x\}$  of $t$ with respect to $x$ (see, for instance,  
Chapter~6 of \cite{Olver} for  a discussion of the Schwarzian derivative and an
illustration of the  role it plays in the analysis of oscillatory differential equations).

\begin{theorem}
\label{theorem:br}
Suppose that $q \in C^\infty\left(\mathbb{R}\right)$ is strictly positive, that
$x(t)$ is defined via (\ref{preliminaries:xt}) and that $p(x)=-2\{t,x\}$.
 Suppose further that there exist positive
constants $\omega$, $\Gamma$ and $\mu$ with
\begin{equation}
\omega > 2 \max\left\{\frac{1}{\mu},\Gamma\right\}
\end{equation}
and
\begin{equation}
\left| \widehat{p}\left(\xi\right) \right|
\leq \Gamma \exp\left(-\mu \left|\xi\right|\right) \ \ \ \mbox{for all}\ \ \
\xi \in \mathbb{R}.
\label{oscillatory:pbound}
\end{equation}
Then there exist smooth functions $\nu$ and $\delta$ such that:
\begin{enumerate}
\item
$\begin{aligned}
\left| \nu(t) \right| \leq
\frac{3\Gamma}{8\mu} \left(1+ \frac{4\Gamma}{\omega}\right)
\exp\left(-\mu \omega\right)
\end{aligned}$
for all $t \in \mathbb{R}$,\\[.1em]

\item
the Fourier transform of $\delta$ is supported on $[-\sqrt{2}\omega,\sqrt{2}\omega]$,\\[.1em]

\item
$\begin{aligned}
\left| \widehat{\delta}(\xi) \right| \leq \left(1+\frac{2\Gamma}{\omega}\right) \frac{\exp\left(-\mu\left|\xi\right|\right)}{4\omega^2-\xi^2}
\end{aligned}$
for all $\left|\xi\right| \leq \sqrt{2}\omega$,\\[.1em]

\item
the function $\alpha$ defined via
\begin{equation}
\alpha(t) = \omega \sqrt{q(t)} \int_0^t \exp\left(\delta(x(s))\right)\, ds
\end{equation}
satisfies the equation
\begin{equation}
\left(\alpha'(t)\right)^2 - \left(\omega^2 q(t) + \nu(t) \right) 
- \frac{3}{4}\left(\frac{\alpha''(t)}{\alpha'(t)}\right)^2 
+ \frac{1}{2} \frac{\alpha'''(t)}{\alpha'(t)} = 0,
\label{preliminaries:kummer}
\end{equation}
\item
the function $r$ defined via
\begin{equation}
r(t) = i\alpha'(t) - \frac{1}{2} \frac{\alpha''(t)}{\alpha'(t)}
\end{equation}
satisfies the equation
\begin{equation}
r'(t) + (r(t))^2 + \left(\omega^2 q(t) + \nu(t) \right) = 0
\ \ \ \mbox{and}
\label{nonoscillatory:perturbed_riccati}
\end{equation}

\item
both 
\begin{equation}
\left\{
\frac{\cos\left(\alpha(t)\right)}{\sqrt{\alpha'(t)}},\ 
\frac{\sin\left(\alpha(t)\right)}{\sqrt{\alpha'(t)}}
\right\}
\end{equation}
and
\begin{equation}
\left\{
\exp\left(\int_0^t r(s)\, ds\right),\ 
\exp\left(\int_0^t \overline{r(s)}\, ds\right)
\right\}
\end{equation}
are bases in the space of solutions of the differential equation
\begin{equation}
y''(t) + \omega^2 \left(q(t)+\nu(t)\right) y(t) = 0,\ \ \  -\infty < t < \infty.
\end{equation}

\end{enumerate}

\end{theorem}

Theorem~\ref{theorem:br} bounds a measure of the  complexity of a solution of the 
perturbed Riccati equation (\ref{nonoscillatory:perturbed_riccati}) in terms of a bound on  the complexity 
of the coefficient $q$.  More explicitly, it states that the Fourier transform of the function
$\delta$  is exponentially decaying with frequency assuming that the same is true
of  $p(t(x))$.   Since $r$ is derived from $q$ and $\delta$, this has the effect
of also bounding the complexity of $r$.  
The perturbation  $\nu$ decays exponentially fast with the frequency
parameter $\omega$, which means that (\ref{nonoscillatory:perturbed_riccati}) becomes indistinguishable
from the Riccati equation for the original differential equation even at very
modest values of $\omega$.  Since  all solutions of the Riccati equation are slowly varying 
when $\omega$ is small, it follows that for all intents and purposes we can always assume the Riccati 
equation has a nonoscillatory solution.  
The bound on the Fourier transform of $\delta$ implies that it can be well-approximated
via a polynomial expansion whose number of terms is independent of $\omega$.
Assuming the same is true of $q$, the solution $r$ of the perturbed Riccati
equation, which is derived from $q$ and $r$, has this property as well.
Finally, we note that while the conditions of Theorem~\ref{theorem:br} appear to be quite stringent,
it can, in fact, be applied to a wide class of differential equations
since it suffices to approximate the coefficient $q$ using one with the required properties.

\end{subsection}

\end{section}

\begin{section}{Analysis of the Discretized Riccati Equation}
\label{section:riccati}

We now prove that when the Riccati equation 
(\ref{introduction:riccati}) is discretized  over an interval $[a,b]$ via a Chebyshev spectral collocation method 
and the Newton-Kantorovich method is applied  to the resulting system of nonlinear algebraic
equations, the iterates converge to a vector which represents a nonoscillatory solution 
of (\ref{introduction:riccati}), provided that the frequency parameter
$\omega$ is sufficiently large.    Throughout this section,  $\left\|\mathbf{v}\right\|$
denotes the  $l^\infty$ norm of the vector $\mathbf{v}$ and $\|T\|$ is the $l^\infty$ 
operator norm of the linear mapping $T: \mathbb{R}^k \to \mathbb{R}^k$.

We begin by discretizing  (\ref{introduction:riccati}) via a standard Chebyshev spectral
collocation method.  
To that end,
we represent the solution $r$ of the Riccati equation by the vector
\begin{equation}
\mathbf{r} = 
\left(\begin{array}{cccc}
r(t_1) &
r(t_2) &
\cdots &
 r(t_k) 
\end{array}
\right)
\end{equation}
of its values at the nodes  $t_1,\ldots,t_k$ of the Chebyshev extremal grid on $[a,b]$, replace the differential
operator by the Chebyshev spectral differentiation matrix
and require that the resulting semi-discrete equation  holds at $t_1,\ldots,t_k$. 
This yields the system of nonlinear algebraic equations
\begin{equation}
\label{riccati:system}
0 = F \left(\mathbf{r}\right)
:= 
\frac{2}{b-a} \Dk \, {\mathbf{r}}
+
\mathbf{r} \circ \mathbf{r}
+ 
\omega^2 \mathbf{q},
\end{equation}
where the vector $q$ is given by
\begin{equation}
\mathbf{q} = 
\left(\begin{array}{cccc}
q(t_1,\omega) &
 q(t_2,\omega) &
\cdots &
 q(t_k,\omega) 
\end{array}
\right).
\end{equation}
and $\mathbf{r}\circ\mathbf{r}$ denotes the Hadmard or entrywise product of the vector
$\mathbf{r}$ with itself.  It can be readily seen that the Fr\'echet derivative of $F$ at $\mathbf{r}$ is 
the matrix
\begin{equation}
F'\left(\mathbf{r}\right)
=
\frac{2}{b-a} \Dk
+
 \mbox{diag}\left(2\mathbf{r}\right).
\end{equation}

We now  apply the Newton-Kantorovich theorem to the system (\ref{introduction:riccati}),
with the initial iterate given by the vector
\begin{equation}
\mathbf{r_0}=
\left(
\begin{array}{ccccc}
\rlg\left(t_1\right) &
\rlg\left(t_2\right)&
\cdots&
\rlg\left(t_k\right)
\end{array}
\right),
\label{riccati:initial}
\end{equation}
where $\rlg$ is defined via
\begin{equation}
\rlg(t) = i \omega \sqrt{q(t,\omega)} - \frac{1}{4} \frac{q'(t)}{q(t)}.
\label{riccati:rlg}
\end{equation}
That is, $\rlg$ is the first order asymptotic approximation of the desired
nonoscillatory solution  of Riccati's equation considered in Section~\ref{section:nonoscillatory}.
We use the notation $\rlg$ because (\ref{riccati:rlg}) is the logarithmic derivative
of the Liouville-Green approximate appearing in (\ref{introduction:lg}).
We could initialize the Newton-Kantorovich iterations with a higher order 
approximation in lieu of $\rlg$; however, it is impractical
to calculate the necessary high-order derivatives of the coefficient $q$ numerically, so we settle
for using the first order approximate as an initial guess.
Under our assumptions on $q$, we have  the bound
\begin{equation}
\omega\, \sqrt{q_{\mbox{\tiny cmin}}} \leq 
\| \mathbf{r_0} \|
\leq \omega + \frac{1}{4 \qmin},
\label{riccati:bound1}
\end{equation}
which holds for all $a \leq t \leq b$ and $\omega \geq \omega_0$. 

As a next step, we establish a  
bound on the   $l^\infty$ operator norm of the  inverse of the Fr\'echet derivative of 
$F$ at the starting point  $\mathbf{r_0}$.  To that end, we observe that
\begin{equation}
\begin{aligned}
\left\| F'\left(\mathbf{r_0}\right)^{-1} \right\|
&=     \left\|  \left(\frac{2}{b-a} \Dk + \mbox{diag}\left(2\mathbf{r_0}\right)\right)^{-1}   \right\| \\
&=
 \left\|\mbox{diag}\left(2\mathbf{r_0}\right)^{-1}  \left( \frac{2}{b-a} \Dk\,\mbox{diag}\left(2\mathbf{r_0}\right)^{-1} + I \right)^{-1}   \right\|\\
&=
 \left\|\mbox{diag}\left(2\mathbf{r_0}\right)^{-1}  \left(I-S \right)^{-1}   \right\|\\
&\leq
 \left\|\mbox{diag}\left(2\mathbf{r_0}\right)^{-1}\right\|
 \left\| I+S +S^2 +S^3 + \cdots   \right\|\\
&\leq
 \left\|\mbox{diag}\left(2\mathbf{r_0}\right)^{-1}\right\|\  \frac{1}{1-\|S\|},
\end{aligned}
\label{riccati:bound3}
\end{equation}
where we have let 
\begin{equation}
S = -\frac{2}{b-a} \Dk\ \mbox{diag}\left(2\mathbf{r_0}\right)^{-1}
\end{equation}
and assumed that the Neumann series for  $(I-S)$ converges.
It is easy to see that this assumption is true when $\omega$ is sufficiently large;
indeed,  (\ref{riccati:bound1}) gives us 
\begin{equation}
\left\| S  \right\| \leq \frac{\|\Dk\|}{\omega (b-a)\sqrt{\qmin}},
\end{equation}
and the right-hand side of the inequality converges to $0$ as $\omega \to \infty$.
So we can choose $\omega$ large enough such that
\begin{equation}
\frac{1}{1-\|S\|} \leq \frac{1}{2},
\label{riccati:bound4}
\end{equation}
which ensures the convergence of the Neumann series and provides us with a bound that  simplifies 
what follows.  Now  (\ref{riccati:bound1}) also implies that
\begin{equation}
\left\| \mbox{diag}\left(2\mathbf{r_0}\right)^{-1}  \right\| \leq \frac{1}{2\omega  \sqrt{\qmin} },
\end{equation}
and combining this with (\ref{riccati:bound4}) yields the desired bound
\begin{equation}
\begin{aligned}
\left\| F'\left(\mathbf{r_0}\right)^{-1}  \right\|
\leq \frac{1}{ \omega  \sqrt{\qmin}} 
\end{aligned}
\label{riccati:bound5}
\end{equation}
on the operator norm of the inverse of the Fr\'echet derivative of $F\left(\mathbf{r_0}\right)$.

We next derive two bounds which will help us to show that fourth and fifth conditions in
Theorem~\ref{theorem:nk} are met.  The assumptions (\ref{introduction:qbound1}) and (\ref{introduction:qbound2})
on $q$ imply that
\begin{equation}
\sup_{a \leq t \leq b} \left|r_{\mbox{\tiny LG}}'(t)  + \left(r_{\mbox{\tiny LG}}(t)\right)^2 + \omega^2 q(t,\omega)\right| 
\leq 
\frac{5}{16} \frac{1}{q_{\mbox{\tiny min}}^2} +  \frac{1}{4} \frac{1}{q_{\mbox{\tiny min}}}.
\label{riccati:rlgbound}
\end{equation}
%
For $k$ sufficiently large, 
\begin{equation}
\left\|\frac{2}{b-a} \mathscr{D}_k\, \mathbf{r_0} -
\left(
\begin{array}{cccccccc}
\rlg'(t_1) & \rlg'(t_2) & \cdots &  \rlg'(t_k)
\end{array}
\right)\right\|
\end{equation}
can be bounded independently of $\omega$.
This together with (\ref{riccati:rlgbound})
implies that there exists a constant $C$ which is independent of $\omega$ and such that 
\begin{equation}
\left\|F\left(\mathbf{r_0}\right) \right\| \leq C.
\label{riccati:fineq}
\end{equation}
Combining (\ref{riccati:bound5}) with (\ref{riccati:fineq})  gives us 
\begin{equation}
\left\| F'\left(\mathbf{r_0}\right)^{-1}  F\left(\mathbf{r_0}\right)\right\|
\leq
\frac{C}{ \omega  \sqrt{\qmin}}.
\label{riccati:bound6}
\end{equation}
Now we observe that 
\begin{equation}
\|F'\left(\mathbf{r}\right) -  F'\left(\mathbf{s}\right)\|
= \left\| 
2 \mbox{diag}\left(2\mathbf{r}\right) 
-2 \mbox{diag}\left(2\mathbf{s}\right) 
\right\|
\leq 2 \|\mathbf{r}-\mathbf{s}\|
\label{riccati:bound6.5}
\end{equation}
for all $\mathbf{r}$ and $\mathbf{s}$ in $\mathbb{R}^k$.  It follows from this and 
(\ref{riccati:bound5}) that 
%
\begin{equation}
\left\| F'\left(\mathbf{r_0}\right)^{-1}\left(
 F'\left(\mathbf{r}\right) -
 F'\left(\mathbf{s}\right)
\right)
\right\|
\leq 
\frac{2}{\omega \sqrt{q_{\mbox{\tiny min}}}}
\left\|\mathbf{r} - \mathbf{s}\right\|
\label{riccati:bound8}
\end{equation}
for all $\mathbf{r}$ and $\mathbf{s}$ in $\mathbb{R}^k$.

Having established the bounds (\ref{riccati:bound6}) and (\ref{riccati:bound8}), 
we are now in a position to choose the constants $\eta$ and $\lambda$ in Theorem~\ref{theorem:nk}
and show that all of its conditions are satisfied.  
We note first that the function $F$ is clearly continuously differentiable
everywhere, so we can take $\Omega=\mathbb{R}^k$ and this implies the second condition
of the theorem is satisfied regardless of our choice of $\eta$.
Earlier, we showed that the  Fr\'echet derivative of $F$ at $\mathbf{r_0}$  is invertible, 
so the first condition is also  satisfied. We now choose the constants to be
\begin{equation}
\eta =  \frac{\omega \sqrt{\qmin}}{2} \ \ \ \mbox{and} \ \ \ 
\lambda =  \frac{C}{ \omega  \sqrt{\qmin}}.
\end{equation}
Then (\ref{riccati:bound6}) is the fourth condition of the theorem, and 
 (\ref{riccati:bound8}) implies the fifth condition.
Moreover, it is clear that the third condition is also satisfied when $\omega$ is sufficiently large.

%
%

Having shown that all of the conditions of Theorem~\ref{theorem:nk}
are satisfied, it now follows that the sequence of iterates $\{\mathbf{r_n}\}$ 
defined by (\ref{riccati:initial}) and 
\begin{equation}
\mathbf{r_{n+1}} = \mathbf{r_n} - F'\left(\mathbf{r_n}\right)^{-1}F\left(\mathbf{r_n}\right)
\end{equation}
converges to a solution $\mathbf{r^*}$ of the discretized Riccati equation.
The theorem also gives us the following bound on the rate of convergence of this
sequence:
\begin{equation}
\left\| 
{\mathbf{r^*}} - {\mathbf{r_n}}
\right\|
\leq
\frac{\eta}{2^n}
\left(
1-\sqrt{1-\frac{2\lambda}{\eta}}
\right)^{2^n}
\leq
\frac{\eta}{2^n} \left(\frac{2\lambda}{\eta}\right)^{2^{n}}
= \mathcal{O}\left(\frac{1}{\omega^{2^{n+1}}}\right)
\ \ \ \mbox{as}\ \omega\to\infty.
\end{equation}

We have only established that 
the iterates  $\mathbf{r_n}$  converge rapidly to a solution of the {\it discretized} 
Riccati equation and, although their limit $\mathbf{r^*}$ corresponds to a polynomial $r^*$ of degree $(k-1)$ 
whose values at the discretization nodes $t_1,\ldots,t_k$ are given by the 
values of the vector $\mathbf{r^*}$, it is not immediately obvious that $r^*$
 approximates a solution of the continuous Riccati equation.
However, from the proof in Section~\ref{section:nonoscillatory1},
we know that the conditions (\ref{introduction:qbound1}) and (\ref{introduction:qbound2})
imply that  there exists a nonoscillatory solution $\rnon$ of the Riccati
equation which can be well-approximated by a polynomial expansion
using a number of terms which is independent of $\omega$, at least
for large enough frequencies.
This means that if the number of collocation points $k$ suffices and the frequency
$\omega$ is large enough, then
\begin{equation}
\mathbf{\rnon} = 
\left(\begin{array}{cccc}
\rnon(t_1) &
\rnon(t_2) &
\cdots &
 \rnon(t_k) 
\end{array}
\right)
\end{equation}
closely approximates a solution of the discretized Riccati equation.
In other words, the nonoscillatory solution of the continuous Riccati equation
yields an approximate solution of the discretized equation.
Moreover, it is clear from (\ref{nonoscillatory:rm}) that 
$\mathbf{\rnon}$ is in the ball $B_\eta(\mathbf{r_0})$ for large enough values of $\omega$.
Since $\mathbf{r^*}$ is the unique
solution of the discretized Riccati equation in that ball,
$\mathbf{r^*}$ must closely coincide with $\mathbf{\rnon}$.
It follows that the vector $\mathbf{r^*}$ must closely approximate the vector 
$\mathbf{\rnon}$ and this, of course, implies that the 
polynomial $r^*$ closely approximate the  nonoscillatory solution $\rnon$ of the 
Riccati equation (\ref{introduction:riccati}).

It appears to be difficult to derive an explicit  lower bound on the value of $\omega$ which 
ensures that the Newton-Kantorovich iterates
converge to a vector which represents  a nonoscillatory solution of the Riccati equation.
It is relatively easy, however, to develop a criterion which works well in practice.
The procedure seems to succeed whenever the frequency is high enough to ensure
that each linearized equation which arises has a unique solution.
The linearization of the discretized Riccati operator around the vector $\mathbf{r}$ is
\begin{equation}
\left( \mbox{diag}(2\mathbf{r})+ \frac{2}{b-a} \mathscr{D}_k \right)\mathbf{h} = -F\left(\mathbf{r}\right),
\label{riccati:linearized1}
\end{equation}
and multiplying both sides of this equation by the inverse of the diagonal matrix yields
\begin{equation}
\left( 
I+S_k\right) \mathbf{h} =
 -\mbox{diag}\left(2\mathbf{r}\right)^{-1}F\left(\mathbf{r}\right),
\label{riccati:linearized2}
\end{equation}
where we have let 
\begin{equation}
S_k=\mbox{diag}\left(2\mathbf{r}\right)^{-1} \frac{2}{b-a} \mathscr{D}_k.
\end{equation}
Because we sample $\rlg$ to initialize the Newton-Kantorovich iterations,
$\mathbf{r} \sim i \omega \sqrt{\mathbf{q}}$, and this gives us the approximate bound
\begin{equation}
\left\|S_k
\right\|
\lesssim  \frac{\left\|\mathscr{D}_k\right\|}{\omega\sqrt{\qmin} (b-a)}.
\label{riccati:sbound}
\end{equation}
We have found (\ref{riccati:sbound}) to be an excellent estimate of the $l^\infty$ operator norm of $S_k$
and when $\|S_k\| <1$, Theorem~\ref{theorem:fp} ensures that 
the linearized equations which arise in the course of the
applying the Newton-Kantorovich method to Riccati's equation 
can be rapidly solved via a fixed point iteration.  

Accordingly, it is tempting to use
\begin{equation}
\frac{\left\|\mathscr{D}_k\right\|}{\omega\sqrt{\qmin} (b-a)} < 1
\label{riccati:criterion0}
\end{equation}
as a criterion for deciding whether the interval $[a,b]$ is in the high-frequency regime
or not.  However, the spectral radius $\rho(S_k)$ of $S_k$ is considerably smaller
than its $l^\infty$ norm, a property inherited from the spectral
differentiation matrix $\mathscr{D}_k$.  Since given 
any $\epsilon > 0$, there exists a matrix norm $\|\cdot\|_\epsilon$ such that 
\begin{equation}
\rho(S_k) \leq \|S_k\|_\epsilon \leq \rho(S_k)+\epsilon,
\end{equation}
it is the spectral radius  of $S_k$ which really controls the behavior
of fixed point iterations for the linearized equation (\ref{riccati:linearized2}).
Consequently, instead of (\ref{riccati:criterion0}), we use the criterion
\begin{equation}
\omega\sqrt{\qmin} (b-a) > \mbox{thresh},
\label{riccati:criterion}
\end{equation}
where $\mbox{thresh}$ is a user-specified threshold parameter which must be adapted
to the choice of $k$.    Assuming $\mbox{thresh}$ is properly set,  (\ref{riccati:criterion})
ensures the existence of a  unique solution
of (\ref{riccati:linearized2}) and  the  rapid convergence of the sequence
$\{\mathbf{h_n}\}$ of fixed point iterates defined via
\begin{equation}
\mathbf{h_0} = 0 \ \ \ \mbox{and}\ \ \ \mathbf{h}_{n+1}  = -\mbox{diag}\left(2 \mathbf{r}\right)^{-1} \frac{2}{b-a} \mathscr{D}_k\,  \mathbf{h_n}
-\mbox{diag}\left(2 \mathbf{r}\right)^{-1} F\left(\mathbf{r}\right)
\label{riccati:fixed}
\end{equation}
to that solution.

\end{section}

\begin{section}{Numerical Algorithm}
\label{section:algorithm}

In this section, we describe our algorithm for the construction of 
of a slowly-varying trigonometric phase function $\alpha$ such that
\begin{equation}
\left\{
\frac{\cos\left(\alpha(t)\right)}{\sqrt{\alpha'(t)}},\ \ 
\frac{\sin\left(\alpha(t)\right)}{\sqrt{\alpha'(t)}}
\right\}
\label{algorithm:basis}
\end{equation}
in basis in the space of solutions of the oscillatory differential equation
(\ref{introduction:ode}).
The algorithm takes the following as inputs:
\begin{enumerate}
\item
 the endpoints $a$ and $b$ of the solution domain $[a,b]$,
\\[-.75em]

\item an external subroutine for evaluating the coefficient $q$ at arbitrary points in $[a,b]$,
\\[-.75em]

\item
an integer parameter $k$ specifying the number of points in the Chebyshev collocation grids used to represent
functions,
\\[-.75em]

\item
a real-valued parameter $\epsilon$ which controls the precision of the calculations and
\\[-.75em]

\item
a real-valued parameter $\mbox{thresh}$ which is used to determine  when an interval
is in the high-frequency regime.
\\[-.75em]

\end{enumerate}
In all of the experiments of this paper, we let $k=16$, $\epsilon=1.0\times 10^{-12}$ and $\mbox{thresh}=10$.
The parameter $\mbox{thresh}$ might need to be adjusted if the value of $k$ is modified.

Our algorithm outputs a collection of discretization intervals $[a_1,b_1],\ [a_2,b_2],\ldots,[a_m,b_m]$ such that
\begin{equation}
a=a_1 < b_1 = a_2 < b_2 = a_3 < b_3 = \cdots = a_m < b_m = b,
\label{algorithm:intervals}
\end{equation}
and the values of $\alpha$, $\alpha'$ and $\alpha''$ at the nodes of the $k$-point Chebyshev extremal grid on each
of these  subintervals.    Using this data, the functions $\alpha$, $\alpha'$ and $\alpha''$ can be easily
evaluated at any point on the interval $[a,b]$.  Moreover, any reasonable initial or boundary
value problem for (\ref{introduction:ode}) can be solved by computing the appropriate
linear combination of the basis functions (\ref{algorithm:basis}).

We note that, as written, our algorithm will fail if there is no discretization interval
in the high-frequency regime.  This is a relatively simple problem to overcome, however,
because in this event, all phase functions for (\ref{introduction:ode}) are
slowly varying throughout the solution domain, and a suitable one
can be constructed by solving the Riccati equation with essentially
arbitrary initial conditions.

The algorithm operates in four stages, each of which is described in detail  below.
In the first stage, we form an initial set of discretization intervals
which suffices to represent the coefficient $q$.
In the second, we traverse traverse this initial collection of discretization
intervals from left-to-right, solving Riccati's equation over each interval in the high-frequency regime
and solving initial value problems for Appell's equation when possible in order to extend the solution
of Riccati equations into the low-frequency regime.  
During the the second stage,  the collection of discretization intervals is adaptively refined  as necessary
in order to represent the phase function accurately.
In the third stage, we traverse the discretization intervals from right-to-left, solving terminal
value problems for Appell's equation in order to extend the solution.  Again, during this stage,
the collection of discretization intervals is adaptively refined as needed.
In the fourth stage, we integrate the obtained solution of Riccati's equation in order to obtain
the desired trigonometric phase function $\alpha$.

In order to make the description of our algorithm simpler,  we first define
several subprocedures.  The most often used of these is the following ``goodness of fit'' procedure for testing
whether or not a function $f$ is well represented by a $(k-1)$-term Chebyshev expansion
on an interval $[c,d]$.  It consists of the following steps:
\begin{enumerate}

\item
Form the vector
\begin{equation}
\mathbf{f} = \left(
f(t_1), \ldots, f(t_k)
\right).
\end{equation}
of values of the function $f$ at the nodes $t_1,\ldots,t_k$ of the Chebyshev extremal grid on the interval $[c,d]$.
\\[-0.75em]

\item
Apply the matrix $\mathscr{C}_k$ to the vector $\mathbf{f}$
in order to compute coefficients $a_0,\ldots,a_{k-1}$ in the Chebyshev expansion
\begin{equation}
\sum_{j=0}^{k-1} a_j T_j\left(\frac{2}{d-c} t  + \frac{d+c}{d-c} \right)
\label{algorithm:expq}
\end{equation}
which agrees with the function $f$ at the nodes $t_1,\ldots,t_k$.
\\[-0.75em]

\item
If the quantity
\begin{equation}
\frac
{ \displaystyle \max \left\{|a_{k-2}|,|a_{k-1}| \right\} } 
{ \displaystyle \max_{j=0,\ldots,k-1} |a_j| }
\end{equation}
is less than the input parameter $\epsilon$, we regard the  expansion
(\ref{algorithm:expq}) as a good approximation of the function $f$.  Otherwise, we regard
it as poor approximation.

\end{enumerate}

The next subprocedure we describe is our method for solving Riccati's equation over an
interval [c,d] in the high-frequency regime.  It is a fairly straightforward application
of the Newton-Kantorovich method with the linearized equations solved 
via the fixed point iteration (\ref{riccati:fixed}).  The only unusual
aspect of our implementation is that we always use the second iterate $\mathbf{h_2}$
in the sequence defined by (\ref{riccati:fixed}) to approximate solutions
of the linearized equations.
This suffices because the fixed point scheme 
converges rapidly and, when applying the Newton-Kantorovich scheme, it
is only necessary to compute the solution of the linearized equation
defining the $(i+1)^{st}$ iterate with accuracy on the order of   
$\| \mathbf{r_{i+1}}-\mathbf{r_{i}} \|$, where $\{\mathbf{r_i}\}$
is the sequence of iterates obtained by solving the linearized equations exactly.
Here are the steps of the procedure in detail:
\begin{enumerate}

\item
Form the vectors
\begin{equation}
\mathbf{q} = 
\left(\begin{array}{cccc}
q(\omega,t_1) &
 q(\omega,t_2) &
\cdots &
 q(\omega,t_k) 
\end{array}
\right)
\end{equation}
and
\begin{equation}
\mathbf{r_0}=
\left(
\begin{array}{ccccc}
\rlg\left(t_1\right) &
\rlg\left(t_2\right)&
\cdots&
\rlg\left(t_k\right)
\end{array}
\right),
\end{equation}
where  $t_1,\ldots,t_k$ are the nodes of the Chebyshev extremal grid on the interval $[c,d]$
and $\rlg$ is defined in (\ref{riccati:rlg}).
Also, set the integer $i$, which is the index of the current Newton-Kantorovich iteration, to $0$.
\\[-.75em]

\item
Compute the residual 
\begin{equation}
F\left(\mathbf{r_i} \right)=  \frac{2}{d-c}\mathscr{D}_k\, \mathbf{r_i} + \mathbf{r_i}\cdot \mathbf{r_i} + \omega^2  \mathbf{q}.
\end{equation}
\\[-.75em]

\item
Approximate the solution of the linearized problem
\begin{equation}
\left(\frac{2}{d-c} \mathscr{D}_k  + \mbox{diag}\left(2 \mathbf{r_i}\right) \right) \mathbf{h} = -F\left(\mathbf{r_i}\right)
\label{algorithm:linearized}
\end{equation}
via
\begin{equation}
\begin{aligned}
\mathbf{h}  &= \mbox{diag}\left(2 \mathbf{r_i}\right)^{-1} 
\left(\mbox{diag}\left(2 \mathbf{r_i}\right)^{-1}  \frac{2}{d-c} \mathscr{D}_k\,   -I \right)F\left(\mathbf{r_i}\right).
\end{aligned}
\end{equation}
This is the second iterate $\mathbf{h_2}$ in the fixed point scheme (\ref{riccati:fixed}).
\\[-.75em]

\item
Form the $(i+1)^{st}$ iterate $\mathbf{r_{i+1}} = \mathbf{r_i} + \mathbf{h}$.
\\[-.75em]

\item
If $\|\mathbf{h}\| > \epsilon \|\mathbf{r_i}\|$ then continue the Newton-Kantorovich iterations by  incrementing $i$ and going to Step~2.
If $\|\mathbf{h}\| \leq \epsilon \|\mathbf{r_i}\|$, then terminate the iterative procedure and regard $\mathbf{r} := \mathbf{r_{i+1}}$ as
the obtained solution of the discretized Riccati equation.
\\[-.75em]
\end{enumerate}

The next subprocedure we describe is an integral equation method for solving an initial value problem
for Appell's equation (\ref{variants:appell}) on the interval $[c,d]$.
When using Chebyshev spectral techniques to solve initial value problems
for differential equations,  we prefer integral formulations to differential
formulations because the former provide a natural way to incorporate the initial
conditions into the problem, while differential formulations require various
ad hoc procedures for enforcing the initial conditions.
Letting
\begin{equation}
m(t) = m(c) + m'(c) (t-c) + m''(c) \frac{(t-c)^2}{2} + \int_c^t \frac{(t-s)^2}{2}\sigma(s)\, ds
\end{equation}
in (\ref{variants:appell}) yields the integral equation
\begin{equation}
\begin{aligned}
&\sigma(t) + 4 \omega^2 q(t,\omega) \int_c^t \int_c^{s_2}  \sigma(s_1)\, ds_1\, ds_2 
+ 2  \omega^2  q'(t) \int_c^t \int_c^{s_3}  \int_c^{s_2} \sigma(s_1)\, ds_1\, ds_2\, ds_3
\\
& \hskip 1em = -4  \omega^2  q(t,\omega) \left( m'(c) + m''(c) (t-c)\right) \\
& \hskip 2em
 -2  \omega^2  q'(t,\omega) \left(m(c) + m'(c) (t-c) + m''(c) \frac{(t-c)^2}{2}\right).
\end{aligned}
\label{algorithm:inteq}
\end{equation}
We discretize and solve it as follows:
\begin{enumerate}

\item
Form the vector
\begin{equation}
\mathbf{q} = 
\left(\begin{array}{cccc}
q(t_1,\omega) &
q(t_2,\omega) &
\cdots &
 q(t_k,\omega) 
\end{array}
\right),
\end{equation}
where  $t_1,\ldots,t_k$ are the nodes of the Chebyshev extremal grid on the interval $[c,d]$. 
\\[-.75em]

\item
 Apply  the
spectral differentiation matrix $2/(d-c) \mathscr{D}_k$ to $\mathbf{q}$ to form the 
vector $\mathbf{qp}$ of values of the derivatives of $q$ at the nodes $t_1,\ldots,t_k$.
\\[-.75em]

\item
Form the $k \times k$ coefficient matrix for the discretized version 
of (\ref{algorithm:inteq}) by letting
\begin{equation}
A = I_k + \mbox{diag}\left(4  \omega^2 \mathbf{q}\right)  \left( \frac{d-c}{2}\mathscr{I}_k\right)^2
        + \mbox{diag}\left(2  \omega^2 \mathbf{qp}\right) \left( \frac{d-c}{2} \mathscr{I}_k \right)^3
\end{equation}
where $I_k$ is the $k\times k$ identity matrix.
\\[-.75em]

\item

Form the vector
\begin{equation}
\begin{aligned}
\mathbf{y} = &-4\omega^2 \left(  m'(c) \mathbf{q}  + m''(c) \mathbf{q}\circ \mathbf{tc} \right)\\
&-2\omega^2 \left( m(c) \mathbf{qp}  + m '(c)  \mathbf{qp} \circ \mathbf{tc}  
+  \frac{1}{2} m''(c) \mathbf{qp}\circ \mathbf{tc} \circ \mathbf{tc}  \right),
\end{aligned}
\end{equation}
where
\begin{equation}
\mathbf{tc} = 
\left(\begin{array}{cccc}
(t_1-c) & (t_2-c) & \cdots & (t_k-c)
\end{array}
\right)
\end{equation}
and $\mathbf{v}\circ\mathbf{w}$ represents pointwise product of the
vectors $\mathbf{v}$ and $\mathbf{w}$. 
\\[-.75em]

\item
Use a standard method to solve the system of linear  equations
$A \pmb{\sigma} = \mathbf{y}$.
\\[-.75em]

\item
Compute vectors $\mathbf{m}$  and $\mathbf{mp}$
representing the solution of Appell's equation and its derivative via the formulas
\begin{equation}
\begin{aligned}
\mathbf{m} &= m(c) \mathbf{1} + m'(c) \mathbf{tc} + \frac{1}{2} m''(c) \mathbf{tc} \circ \mathbf{tc}  + \left(\frac{d-c}{2} \mathscr{I}_k\right)^3 \pmb{\sigma}
\ \ \ \mbox{and}\\\
\mathbf{mp} &= m'(c) \mathbf{1} + m''(c) \mathbf{tc}   + \left(\frac{d-c}{2} \mathscr{I}_k\right)^2 \pmb{\sigma},\\
\end{aligned}
\end{equation}
where $\mathbf{1}$ is the vector whose entries are all $1$'s and $\mathbf{tc}$
is as before.

\end{enumerate}
We use a completely analogous procedure to solve terminal boundary value problems
for Appell's equation.  Because of the great similarities between
that procedure and the one just described, we omit a detailed description of it.

We are now in a position to describe each stage of our algorithm.

\begin{subsection}{Stage one: adaptive discretization of the coefficient $q$}

In this stage, we construct an initial list of intervals $[a_1,b_1],\ [a_2,b_2], \ldots, [a_m,b_m]$
which suffice to discretize the coefficient $q$.  
In addition to this list of discretization intervals,
a list of intervals to process is maintained.
In the first instance,  the list of discretization intervals
is empty and the list of intervals to process contains only $[a,b]$.
This stage proceeds by executing the following steps until the list
of intervals  to process is empty:

\begin{enumerate}
\item
Remove an interval $[c,d]$ form the list of intervals to process.
\\[-.75em]

\item
Check the accuracy with which the $q$ is represented
via a $k$-term Chebyshev expansion over $[c,d]$ using
the goodness of fit procedure.  If it is deemed to be well represented,
put $[c,d]$ into the list of discretization intervals.
Otherwise,
put the intervals $[c,(c+d)/2]$ and $[(c+d)/2,d]$ into the
list of intervals to process.    

\end{enumerate}

\end{subsection}

\begin{subsection}{Stage two: solution of the Riccati equation and left-to-right sweep}

In this stage, we traverse the list $[a_1,b_1],\ldots,[a_m,b_m]$ of discretization intervals
formed in the previous stage from left-to-right.  We solve Riccati's equation on any interval
in the high-frequency regime, and solve an initial value problem for Appell's equation
in order to extend the phase function whenever possible.
During this stage, we adaptively refine the list of discretization intervals as needed 
to ensure the phase function is accurately represented.  
To this end, we maintain two lists, one containing the intervals 
which need to be  processed and the other specifying the new collection of discretization 
intervals created during this stage.  The list of intervals to process is initialized
with all of the discretization intervals constructed in the preceding phase,
and the following sequence of steps is executed
until the list of intervals to process is empty:

\begin{enumerate}

\item
Remove from the list of intervals to process the left-most interval
$[c,d]$.
\\[-.75em]

\item
Form the vector $\mathbf{q}$ of the values of the coefficient at the
nodes $t_1,\ldots,t_k$ of the Chebyshev extremal grid on $[c,d]$
and let $\qmin$ denote the minimum value of $q$ at those nodes.
\\[-.75em]

\item
Compute the quantity $\gamma = \omega\sqrt{\qmin} (b-a)$.
\\[-.75em]

\item
If $\gamma > \mbox{thresh}$, execute the following sequence of steps:
\\[-.75em]


\begin{enumerate}

\item
Compute a solution $\mathbf{r}$ of the Riccati
equation using the subprocedure described above.
\\[-.5em]

\item
Let $\pmb{ap}$ be the real-valued vector whose entries
are the imaginary parts of the entries of $\mathbf{r}$, and 
let $\pmb{app}$ be the pointwise produce of the vector
$-2 \pmb{ap}$ and the real part of $\mathbf{r}$.  
According to Formula~(\ref{variants:ralpha}), these are the computed values of the derivative $\alpha'$ and second derivative
$\alpha''$ of the desired trigonometric phase function $\alpha$.
\\[-.5em]

\item
Perform the goodness of fit procedure on $\pmb{ap}$.  If it is 
well-approximated by $(k-1)$-term Chebyshev expansion, then
add the interval $[c,d]$ to the list of output intervals.
Otherwise, add the intervals $[c,(c+d)/2]$ and
$[(c+d)/2,d]$ to the list of intervals
to process.  In either case, goto Step~1
of the procedure of this stage.  
\\[-.5em]
\end{enumerate}

\item
If $\gamma \geq \mbox{thresh}$ and the functions $\alpha'$ and $\alpha''$
have already been computed over  the interval immediately
to the left of $[c,d]$, then we solve an initial value problem for 
Appell's equation to construct 
the vectors $\mathbf{ap}$ and $\mathbf{app}$ of values
of $\alpha'$ and $\alpha''$ on $[c,d]$ using the following sequence of steps:
\\[-.75em]


\begin{enumerate}

\item
Let  $apval$ and $appval$ denote the values of $\alpha'(c)$  and $\alpha''(c)$ (these are already
known since $c$ is the right-hand endpoint of the interval to the left of $[c,d]$)  and
let  $qval$ be the value of the coefficient $q$ at the point $c$.
Use Kummer's equation (\ref{introduction:kummer}) 
to form an estimate apppval of the value of $\alpha'''(c)$ via the formula
\begin{equation}
apppval = \frac{ 4\, qval \times apval^2-4apval^4+3\, appval^2 } { 2\, apval}.
\end{equation}
\\[-.5em]

\item
Calculate the initial values of the solution $m$ of 
Appell's equation such that  $1/m$ extends 
the derivative $\alpha'$ of the trigonometric phase function from
the preceding interval to the current one.
More explicitly, let 
\begin{equation}
\begin{aligned}
m(c)   &= \frac{1}{apval},\ \ \ 
m'(c)  &= -\frac{appval}{apval^2}\ \ \ \mbox{and} \ \ \ 
m''(c) &= 2\frac{appval^2}{apval^3} - \frac{apppval}{apval^2}.
\end{aligned}
\end{equation}
\\[-.5em]

\item
Solve the initial value problem for Appell's equation using the procedure described above.
This results in vectors $\mathbf{m}$  and $\mathbf{mp}$ which give
the values of the solution of Appell's equation $m(t)$ and its derivative $m'(t)$
at the Chebyshev nodes $t_1,\ldots,t_k$.
\\[-.5em]

\item
Compute the vectors $\mathbf{ap}$ and $\mathbf{app}$, whose entries give the 
values of $\alpha'$ and $\alpha''$ at the Chebyshev nodes $t_1,\ldots,t_k$ on $[c,d]$,
using vectors $\mathbf{m}$ and $\mathbf{mp}$.  To be more explicit, 
$\mathbf{ap}$ is the vector whose entries are the reciprocals
of those of $\mathbf{m}$, and 
\begin{equation}
\mathbf{app} = - \mathbf{ap} \circ \mathbf{ap} \circ \mathbf{mp}.
\end{equation}
\\[-.5em]

\item
Perform the goodness of fit procedure on $\pmb{ap}$.  If it is 
well-approximated by $(k-1)$-term Chebyshev expansion, then
add the interval $[c,d]$ to the list of output intervals.
Otherwise, add the intervals $[c,(c+d)/2]$ and
$[(c+d)/2,d]$ to the list of intervals
to process.  In either case, goto Step~1
of the procedure of this stage.  

\end{enumerate}
\vskip 1em


\end{enumerate}

\end{subsection}

\begin{subsection}{Stage three: right-to-left sweep}

At the conclusion of the second state of the procedure, 
we have a refined list $[a_1,b_1],\ldots,[a_m,b_m]$ of discretization intervals
and the vectors $\mathbf{ap}$ and $\mathbf{app}$ have been constructed
over all intervals which lie to the right of an interval in the high-frequency regime.
In this stage, we sweep from right-to-left, solving terminal value problems
for Appell's equation over any interval for which $\mathbf{ap}$ and
$\mathbf{app}$ have not yet been constructed.

As before, we maintain two lists, one containing the intervals 
which need to be  processed and the other specifying the new collection of discretization 
intervals created during this stage.  The list of intervals to process is initialized
with all of the discretization intervals constructed in the preceding phase,
and the following sequence of steps is executed
until the list of intervals to process is empty:

\begin{enumerate}

\item
Remove from the list of discretization intervals to process the right-most interval
$[c,d]$ for which the vectors $\mathbf{ap}$ and $\mathbf{app}$ have not yet been constructed.
\\[-.75em]

\item
Solve a terminal value problem for Appell's equation in order to construct
the vectors $\mathbf{ap}$ and $\mathbf{app}$
We omit the details because the procedure is entirely analogous to that used
in the second stage of the procedure.
\\[-.75em]

\item
Perform the goodness of fit procedure on $\pmb{ap}$.  If it is 
well-approximated by $(k-1)$-term Chebyshev expansion, then
add the interval $[c,d]$ to the list of output intervals.
Otherwise, add the intervals $[c,(c+d)/2]$ and
$[(c+d)/2,d]$ to the list of intervals
to process.  In either case, goto Step~1
of the procedure of this stage.  

\end{enumerate}

\end{subsection}

\begin{subsection}{Stage four: spectral integration}

At the conclusion of the third stage, we have 
the final list $[a_1,b_1],\ldots,[a_m,b_m]$ of discretization
intervals and the vectors $\mathbf{ap}$ and $\mathbf{app}$
have been constructed over each interval.
In this final stage, we use spectral integration
to compute the values of the desired nonoscillatory trigonometric phase function
$\alpha$ at the Chebyshev nodes of each discretization interval.
More explicitly, we initial set $aval=0$ and then, for each $j=1,\ldots,m$ we execute the following steps:
\begin{enumerate}

\item
Compute the vector $\mathbf{a}$ giving the values of $\alpha$ at the Chebyshev nodes on $[a_j,b_j]$ via
the formula

\begin{equation}
\mathbf{a} = aval + \frac{b_j-a_j}{2} \mathscr{I}_k \mathbf{ap}
\end{equation}
\\[-.75em]

\item
We let aval be equal to the last entry of the vector $\mathbf{a}$.

\end{enumerate}

\end{subsection}

\end{section}

\begin{section}{Numerical Experiments}
\label{section:experiments}

In this section, we present the results of numerical experiments which were conducted to illustrate 
the properties of the method of this paper and to compare it with other solvers for oscillatory
differential equations. We implemented our algorithm in Fortran and compiled our code with  version 14.2.1 of the 
GNU Fortran compiler.    All experiments were performed on a desktop computer equipped
with an AMD 9950X processor and 64GB of RAM.  This processor  has 16 cores, but
only one core was utilized in our experiments.

In all of our experiments,  the value of the parameter $k$, which determines the order of the 
Chebyshev expansions used to represent phase functions, was taken to be $16$,
the parameter $\epsilon$  that controls the accuracy of the obtained phase functions
was set to $10^{-12}$ and the parameter $\mbox{thres}$ was set to $10$.

In almost all of our experiments, we tested the accuracy of our method and those we compare it to 
by using it to calculate solutions to initial and boundary value problems for second
order equations. The experiment of   Section~\ref{section:experiments3},
in which the accuracy of the phase functions produced by our algorithm was measured
directly, is the sole exception.
Because the condition numbers of initial and boundary value problems for equations of the form
(\ref{introduction:ode}) grow with $\omega$, the accuracy of any numerical
method used to solve them is expected to deteriorate with increasing  frequency.    
In the case of our algorithm, the phase functions themselves are calculated to high precision, 
but the magnitudes of the phase functions  increase with frequency and accuracy is lost when 
the phase functions  are exponentiated.
One implication is that calculations which involve only the phase functions
and not the solutions of the scalar equation can be performed to high accuracy.
The article \cite{BremerZeros}, for example,  describes
a phase function method  for rapidly computing the zeros of solutions of second order linear ordinary
differential equations to extremely high accuracy.    

To account for the vagaries of modern computing environments, all reported times were obtained 
by averaging the cost of each calculation over 100 runs.

%
%

\begin{subsection}{Comparison with the Modified Magnus expansion method}
\label{section:experiments1}

In this first experiment, we compared the performance of the modified Magnus method of \cite{ModMagnus}, 
which is one of the most widely used methods for solving oscillatory differential equations,
with the algorithm of this paper.  The former is a step method
that combines an exponential integrator with  preconditioning via the solutions of a 
constant coefficient equation obtained by freezing the coefficient matrix of the oscillatory 
equation. We used a fourth order exponential integrator and 
equispaced step sizes in our implementation of the  modified Magnus method,
which was written in Fortran and closely follows the
description provided in \cite{ModMagnus} and \cite{Iserles2}.

For each $n=2^6, 2^7, 2^8, \ldots, 2^{20}$, we used both methods to evaluate the solution
\begin{equation}
L_n(x) = P_n(x) + i \frac{2}{\pi} Q_n(x)
\label{experiments:lnu}
\end{equation}
of Legendre's differential equation
\begin{equation}
(1-t^2) y''(t)  -2t y'(t) + n(n+1) y(t) = 0 
\label{experiments:legendre}
\end{equation}
at a collection of points on the interval $[0,0.9]$.  Because Equation~(\ref{experiments:lnu})
is not in the normal form (\ref{introduction:ode}), we actually applied the algorithm 
of this paper to the differential equation
\begin{equation}
y''(t) + \left(\frac{1}{(1-t^2)^2} + \frac{n(n+1)}{1-t^2} \right) y(t) = 0,
\label{experiments:legnormal}
\end{equation}
which has  $L_n(t)\sqrt{1-t^2}$ as a solution.
We choose to evaluate the Legendre function $L_n$ because its logarithmic derivative
is nonoscillatory.  We discuss this further  in Section~\ref{section:experiments3},
and it can be seen from the graph of the accuracy predicted by its
condition number of evaluation in Figure~\ref{experiments1:figure1}.

For each value of $n$ considered,  we executed the modified Magnus method twice, once
with the step size $h$ taken to be a constant multiple of $n^{-1/2}$ 
and again with the step size taken to be the same constant multiple of $n^{-3/4}$.
The constant was set so as to ensure  $10^{-10}$ relative accuracy for 
the smallest value of $n$ considered.
We struggled to obtain higher accuracy with the modified Magnus
method at high frequencies because it becomes numerically unstable when the step size is 
extremely small. In the case of our algorithm, for each value of $n$ considered,
we evaluated $L_n$ at a collection of 1,000 equispaced points on the 
interval $[0,0.9]$, including $0$ and $0.9$, and recorded the largest relative error which was observed.
For the modified Magnus method, we recorded the largest
relative error which was observed at the steps taken by the solver.
Figure~\ref{experiments1:figure1} gives the results.  The plot on the left
shows   the time taken by each method as a function of $n$, while that on the right 
gives the largest observed relative errors, again as a function of $n$.  
The plot on the right also shows the relative accuracy predicted by the condition 
number of evaluation of the function $L_n$.  
More explicitly, it displays a graph of 
\begin{equation}
\kappa(n) = \max_{j=1,\ldots,1000} \epsilon_0 \left|\frac{t_j L_n'(t_j) }{L_n(t_j)}\right|,
\end{equation}
where $t_1,\ldots,t_{1000}$ are the equispaced nodes on $[0,0.9]$ at which
we evaluated $L_n$ using the algorithm of this paper and 
 $\epsilon_0 \approx 2.220446049250313\times 10^{-16}$
is machine zero.

From Figure~\ref{experiments1:figure1}, we see that the modified Magnus method
maintains roughly constant accuracy when the number of steps scales
as $n^{-3/4}$, but loses accuracy when the step size scales as $n^{-1/2}$.
This is consistent with the theoretical estimates presented in \cite{ModMagnus}.
We note that the method of this paper is several orders of magnitude faster than the modified
Magnus method, even at low frequencies, and it is more than 4 orders of magnitude
faster at the highest frequency considered.

\begin{figure}[h!!!!!!!!!!!!!!!!!!!!]
\hfil
\includegraphics[width=.49\textwidth]{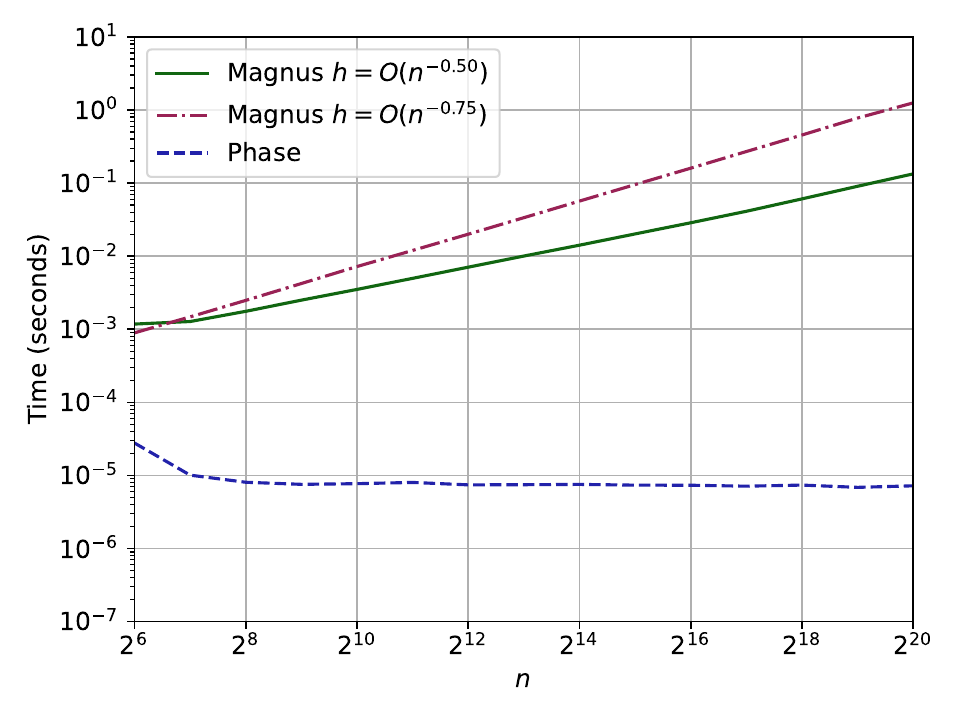}
\hfil
\includegraphics[width=.49\textwidth]{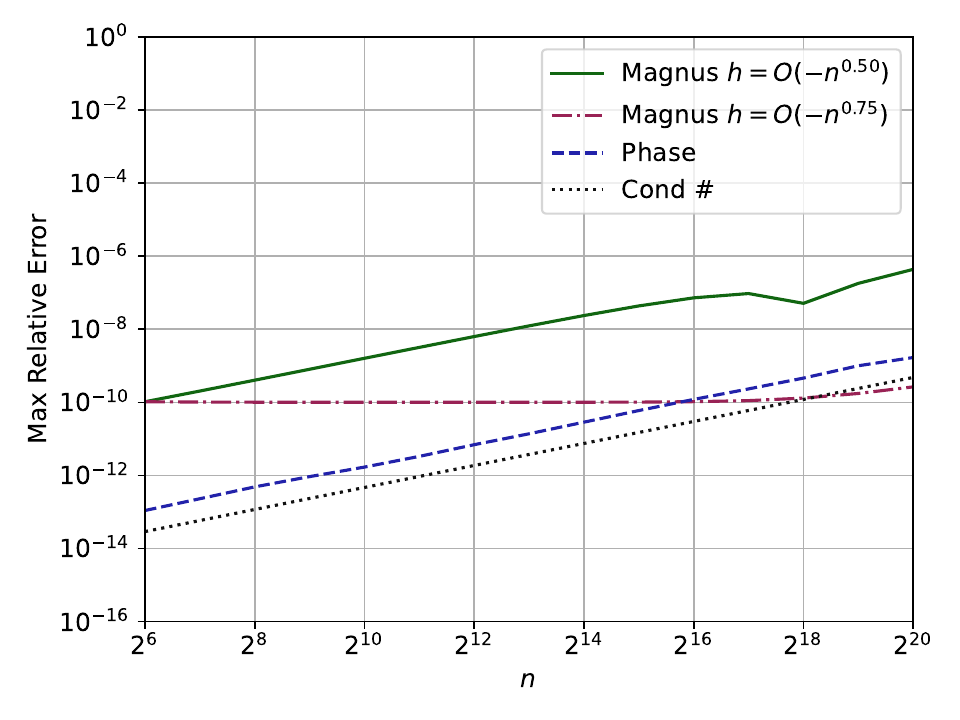}
\hfil

\caption{The results of the experiment of Section~\ref{section:experiments1}
in which the modified Magnus method of \cite{ModMagnus} and the algorithm
of this paper were used to evaluate the Legendre function $L_n(t)$ defined in (\ref{experiments:lnu})
on the interval $[0,0.9]$.  
The plot on the left gives the time in seconds required by each method as a function
of the degree $n$ of the Legendre function, while the plot on the right
gives the maximum relative error observed in the course of evaluating  the Legendre function
using each method. The plot on the right  also shows the accuracy predicted by the 
condition number of evaluation of the Legendre function.
We give the results for the modified Magnus method when the step size $h$
was taken to be a constant multiple of  $n^{-1/2}$, and when it was taken to be a constant multiple of
 $n^{-3/4}$.  
}
\label{experiments1:figure1}
\end{figure}

\end{subsection}

%
%

\begin{subsection}{Comparison with two frequency-independent solvers}
\label{section:experiments2}

In this experiment, we compared the algorithm of this paper with the ARDC
method of \cite{ARDC} and the frequency-independent solver \cite{BremerPhase}
previously developed by one of this paper's authors
by once again evaluating the solution $L_n$ of Legendre's differential
equation (\ref{experiments:legendre}).
More explicitly, for each $n=2^6, 2^7, 2^8, \ldots, 2^{20}$, we used each
method to compute 
%
the Legendre function $L_n$ defined in (\ref{experiments:lnu})
at $100$ equispaced nodes on $[0,0.999]$, including 
the points $0$ and $0.999$.
We evaluated $L_n$ at a relatively small number of nodes 
since the method of \cite{ARDC} is a step scheme which
only returns the values of the solution at a sparse set of 
nodes sampled well below the Nyquist frequency.  This means that the solution
returned by the ARDC method  cannot be accurately interpolated at  arbitrary
points on the solution domain.  In order to evaluate it at 
the  $100$ equispaced nodes we considered, we had to specify them as inputs
to ARDC and doing so increased the cost of the algorithm.  Adding more points would
increase the cost further, and we regard $100$ points as  sufficient
to measure the accuracy of the solution.    Both the algorithm
of this paper and that of \cite{BremerPhase} allow for the solution
to be evaluated at any point on the solution domain.

Figure~\ref{experiments2:figure1} gives the results of this experiment.
The plot on the left gives the time required by each method
as a function of $n$, while that on the right
shows the maximum relative accuracy achieved by each method as a function of $n$, 
as well as the relative accuracy predicted by the condition number of evaluation of $L_n$.
We note that the three methods achieve very similar levels of accuracy
for large values of $n$, while the the ARDC method loses some accuracy
at small values of $n$.  ARDC is noticeably slower than \cite{BremerPhase}, and it
is several orders of magnitude slower than the method of this paper.

\begin{figure}[h!!!!!!!!!!!!!!!!!!!!]
\hfil
\includegraphics[width=.49\textwidth]{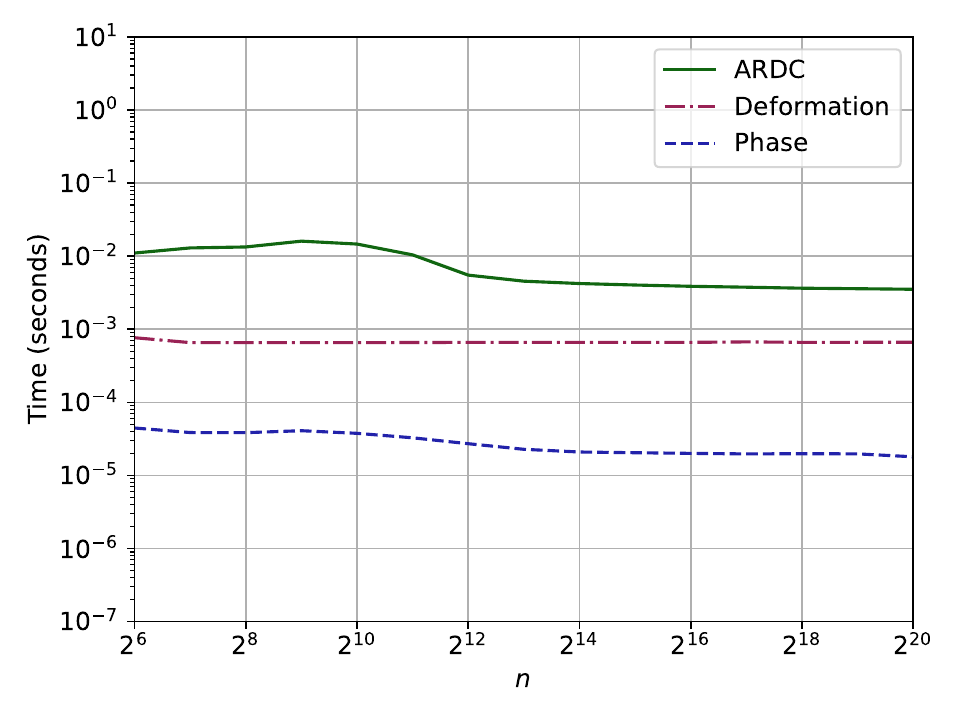}
\hfil
\includegraphics[width=.49\textwidth]{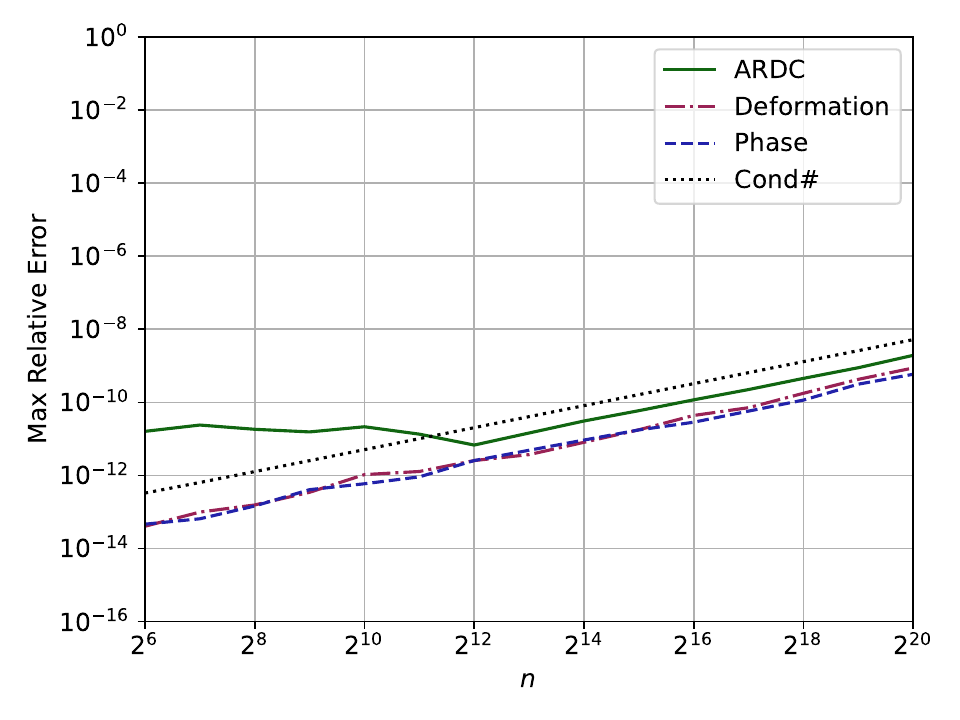}
\hfil

\caption{The results of the experiment of Subsection~\ref{section:experiments2}
in which the method of this paper, the smooth deformation scheme of \cite{BremerPhase}
and  the ARDC method of \cite{ARDC} were used to evaluate the Legendre function
$L_n$ defined in (\ref{experiments:lnu})
on the interval $[0,0.999]$.
The plot on the left gives the time required by each algorithm as a function of the degree
$n$ of the Legendre function,
while the plot on the right shows the maximum relative error observed
when $L_n$ was evaluated
at $100$ equispaced points on the interval $[0,.999]$ using each approach.
The plot on the right also gives the maximum relative error 
predicted by the condition number of evaluation of this function.
}
\label{experiments2:figure1}
\end{figure}

\end{subsection}

%
%

\begin{subsection}{Accuracy of the computed phase functions}
\label{section:experiments3}

In this experiment, we  measured the accuracy with which our 
method  computes phase functions.  For each $n=2^7,2^8,\ldots,2^{21}$, we 
used our algorithm   to construct the derivative of a nonoscillatory
trigonometric phase function for the  normal form (\ref{experiments:legnormal}) of
Legendre's differential equation on the interval $\left[0,1-10^{-7}\right]$.

There is an explicit expression for  the desired solution of Kummer's equation in this case.
Indeed, it can be shown using the formula
\begin{equation}
 \left(P_n(t)\right)^2 + \frac{4}{\pi^2} \left(Q_n(t)\right)^2
= \frac{4}{\pi^2} \int_1^\infty Q_n\left(t^2 + (1-t^2)s\right) \frac{ds}{\sqrt{s^2-1}},
\ \ \ -1 < t < 1,
\label{experiments3:legemod}
\end{equation}
which appears in \cite{durand75}, that the function
on the  left-hand side of (\ref{experiments3:legemod}) is absolutely monotone on $(0,1)$.
Since the functions
\begin{equation}
\sqrt{\frac{\pi}{2}} P_n(t)\sqrt{1-t^2}, \ \ \ \frac{2}{\pi}  Q_n(t)\sqrt{1-t^2}
\label{experiments3:pq}
\end{equation}
are a pair of solutions of (\ref{experiments:legnormal})
whose Wronskian is $1$ (see, for instance, Chapter~3 of \cite{HTFI}), the function
\begin{equation}
\alpha'(t) = \frac{1}{\left(1-t^2\right)\left(\frac{\pi}{2}P_n^2(t) + \frac{2}{\pi} Q_n^2(t) \right)}.
\label{experiments3:kummer}
\end{equation}
is the derivative of the desired nonoscillatory trigonometric phase function
for the normal form of Legendre's differential equation.

For each value of $n$ considered in this experiment,
we measured the accuracy of the  obtained solution of Kummer's equation
by comparing it with the reference solution (\ref{experiments3:kummer})
at  at 1,000 equispaced points on the interval $\left[0,1-10^{-7}\right]$,
including the points $0$ and $1-10^{-7}$.
Figure~\ref{experiments3:figure1} 
gives the largest observed relative error in the computed value of the phase function as a function of $n$,
as well as the time required to construct the phase function.
From these results, we see that the phase function is always computed with
relative accuracy which is greater than the requested
precision $1.0 \times 10^{-12}$.  

\begin{figure}[h!!!!!!!!!!!!!!!!!!!!]
\hfil
\includegraphics[width=.49\textwidth]{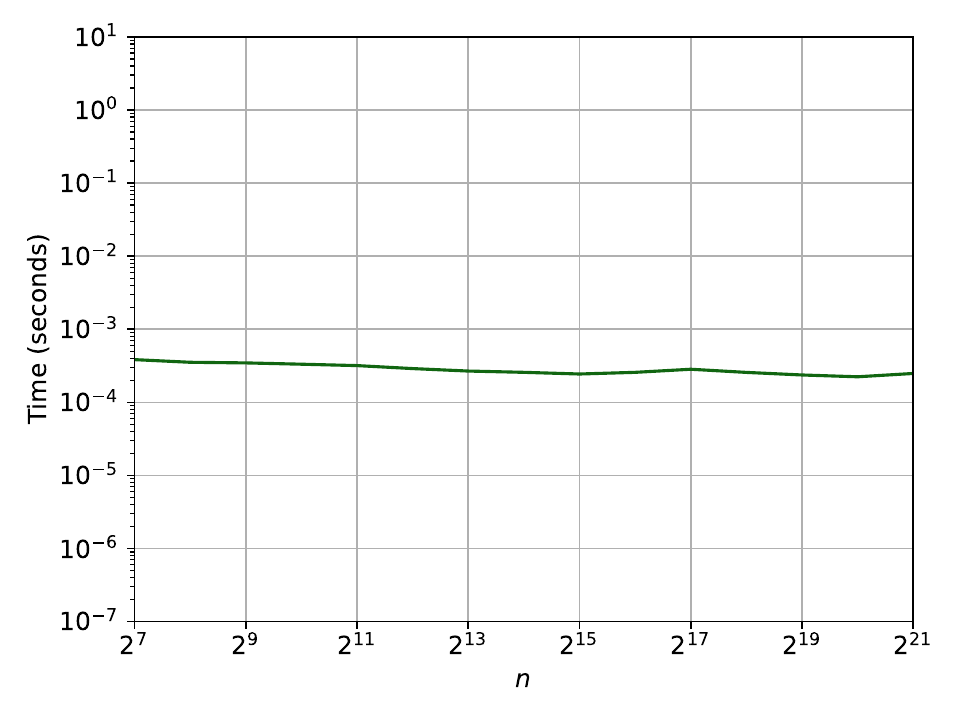}
\hfil
\includegraphics[width=.49\textwidth]{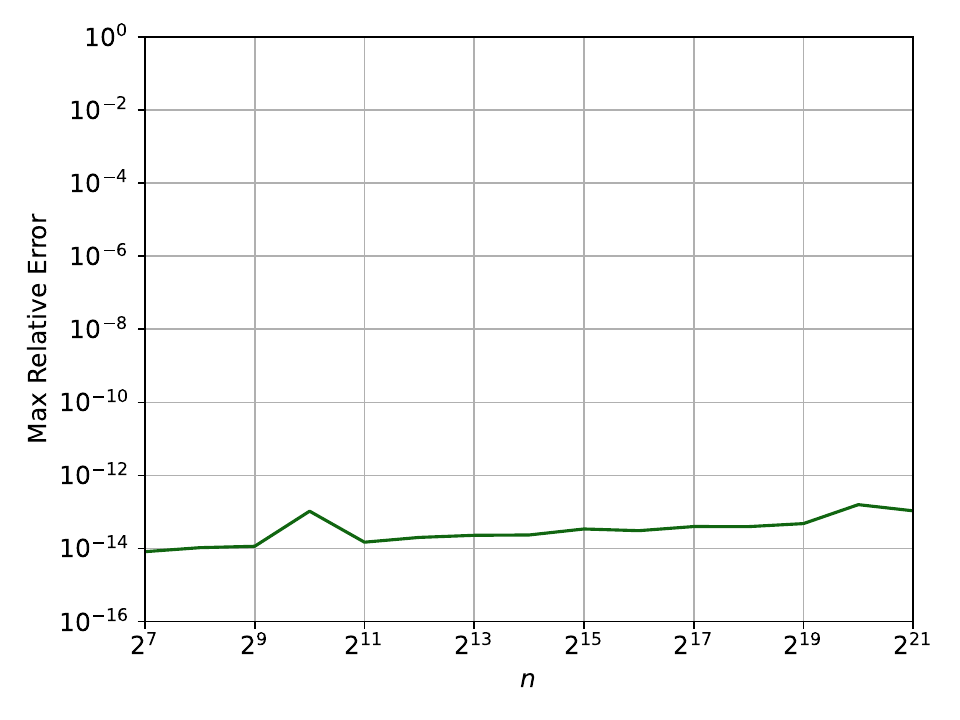}
\hfil

\caption{The results of the experiment of Subsection~\ref{section:experiments3}
in which the accuracy of the phase functions computed by our algorithm
was measured.  The plot on the left 
gives the time taken by our algorithm as a function of the frequency parameter $n$.
The plot on the right gives the maximum relative error
observed while evaluating the solution of Kummer's equation calculated by our algorithm
at 1,000 points in the interval $\left[0,1-10^{-7}\right]$, including
the points $0$ and $1-10^{-7}$.
}
\label{experiments3:figure1}
\end{figure}

\end{subsection}

%
%

\begin{subsection}{Evaluation of Gegenbauer polynomials}
\label{section:experiments4}
In the experiment of this section, we used the scheme of this paper 
to evaluate the Gegenbauer polynomial $C_n^{\alpha}(t)$, which is
the unique  solution of the differential equation
\begin{equation}
(1-t^2) y''(t) - 2(\alpha+1) y(t) + n(n+2 \alpha)y(t) = 0
\label{experiments4:ode}
\end{equation}
that is continuous on $[-1,1]$ and such that 
\begin{equation}
C_n^\alpha(1) = \frac{\Gamma(2\alpha+1)}{\Gamma(2\alpha)\Gamma(n+1)}.
\end{equation}
For each $n=2^6,2^7,2^8,\ldots,2^{20}$ and $\alpha=-0.499,0.25,1.00$,
we computed a phase function for the normal form
\begin{equation}
y''(t) + \left(
\frac{\alpha-\alpha^2+\frac{3}{4}}{(1-t^2)^2}
+
\frac{\left(n+\alpha-\frac{1}{2}\right)\left(n+\alpha+\frac{1}{2}\right)}{1-t^2}
\right) y(t) = 0
\end{equation}
of (\ref{experiments4:ode}) satisfied by $C_n^{\alpha}(t)\left(1-t^2\right)^{(2\alpha+1)/4}$
via the algorithm of this paper and used
it to evaluate $C_n^{\alpha}(x)$ at 1,000 equispaced points on the
interval $(0,0.999)$.  
We compared the obtained values to those computed using the well-known three-term recurrence
relation satisfied by the Gegenbauer polynomials.
The results are shown in Figure~\ref{experiments4:figure1}.  There, 
we give the largest observed absolute error as a function of $n$,
as well as the time required to construct each  phase function.
We measured absolute error rather than relative error in these experiments because $C_n^\alpha$ has zeros
on the interval $(0,1)$.  

\begin{figure}[h!!!!!!!!!!!!!!!!!!!!]
\hfil 
\includegraphics[width=.49\textwidth]{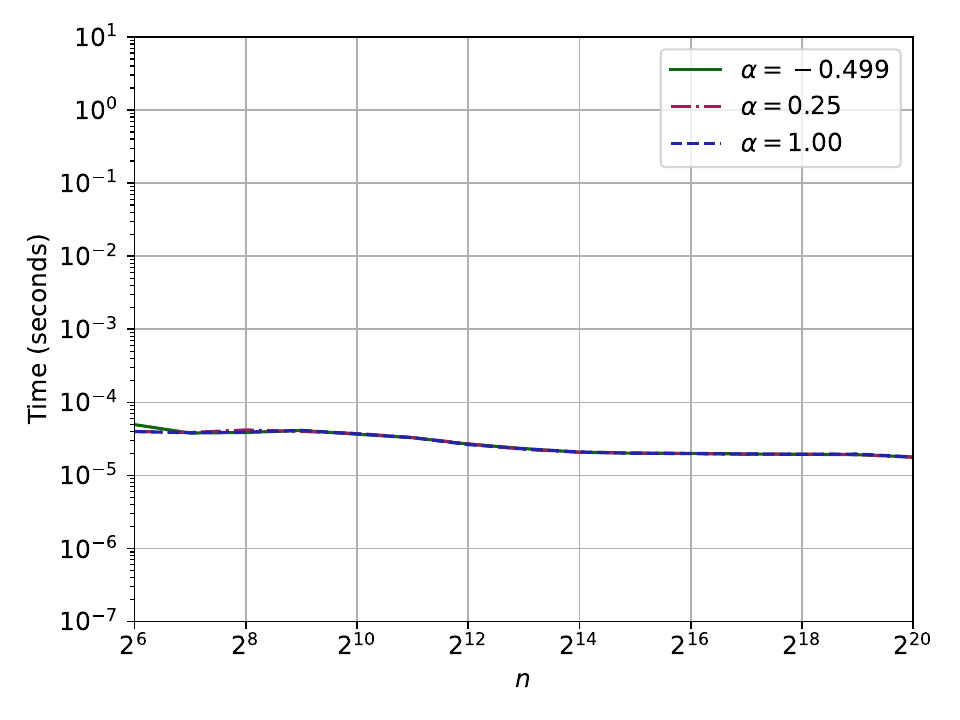}
\hfil
\includegraphics[width=.49\textwidth]{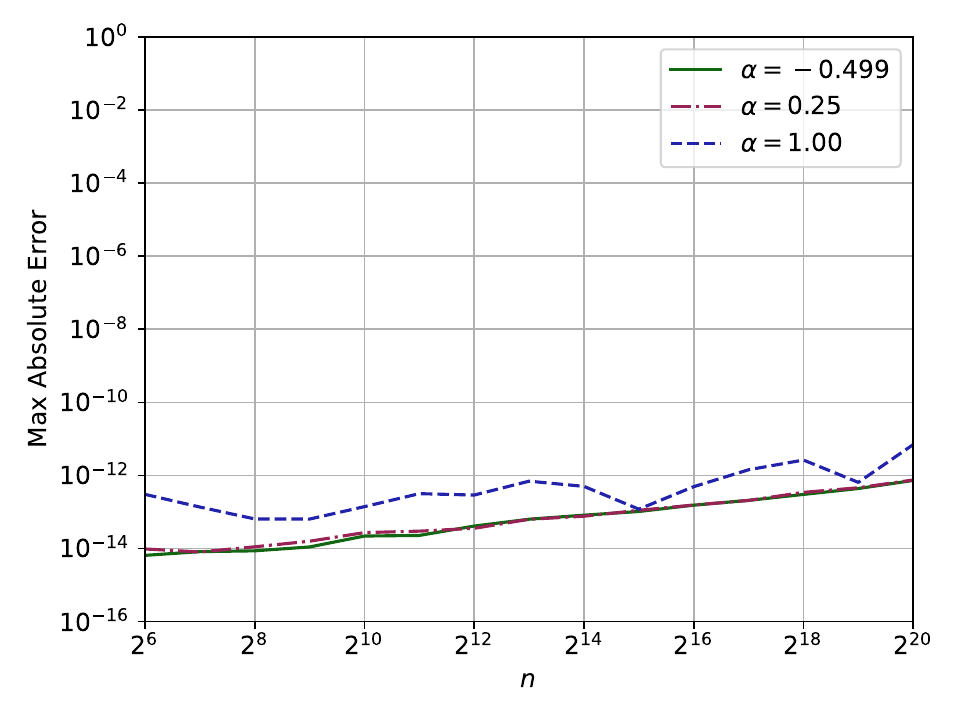}
\hfil

\caption{The results of the experiment of Subsection~\ref{section:experiments4}
in which our algorithm was used to evaluate Gegenbauer polynomials
of orders $\alpha=-0.499,0.25,1.0$.
The plot on the left  gives the time required to construct each phase function.
The plot on the right gives the maximum absolute accuracy observed
while evaluating the Gegenbauer polynomials at 1,000 equispaced
points on the interval $[0,0.999]$.
}
\label{experiments4:figure1}
\end{figure}

\end{subsection}

%
%

\begin{subsection}{A boundary value problem}
\label{section:experiments5}

In this final experiment, 
for each $\omega=2^6,2^7,\ldots,2^{20}$, 
we constructed a phase function for the differential equation
\begin{equation}
y''(t) + 
\omega ^2 \left(\frac{3 t^2
   \omega ^2+t^2 \omega
   +1}{-\left(t^2+1\right)
   \omega +\omega ^2+1}+\frac{2
   e^{-t}}{t^2+\frac{1}{10}}\right)
y''(t) = 0,\ \ \ -1 < t < 1,\ \ \ 
\label{experiments:bvp}
\end{equation}
and used it to evaluate the solution $y$ which satisfied the boundary conditions
 $y(-1) = y(1) = 1$ at 1,000 equispaced points on $[-1,1]$.  We then
applied  a standard adaptive Chebyshev spectral method to these boundary value problems
in order to construct reference solutions.
Figure~\ref{experiments5:figure1} gives the results.  The time required by our algorithm
and by the standard solver are plotted as functions of the frequency $\omega$
on the left, while the maximum observed
absolute error is plotted as a function of $\omega$ on the right.

\begin{figure}[h!!!!!!!!!!!!!!!!!!!!]
\hfil 
\includegraphics[width=.49\textwidth]{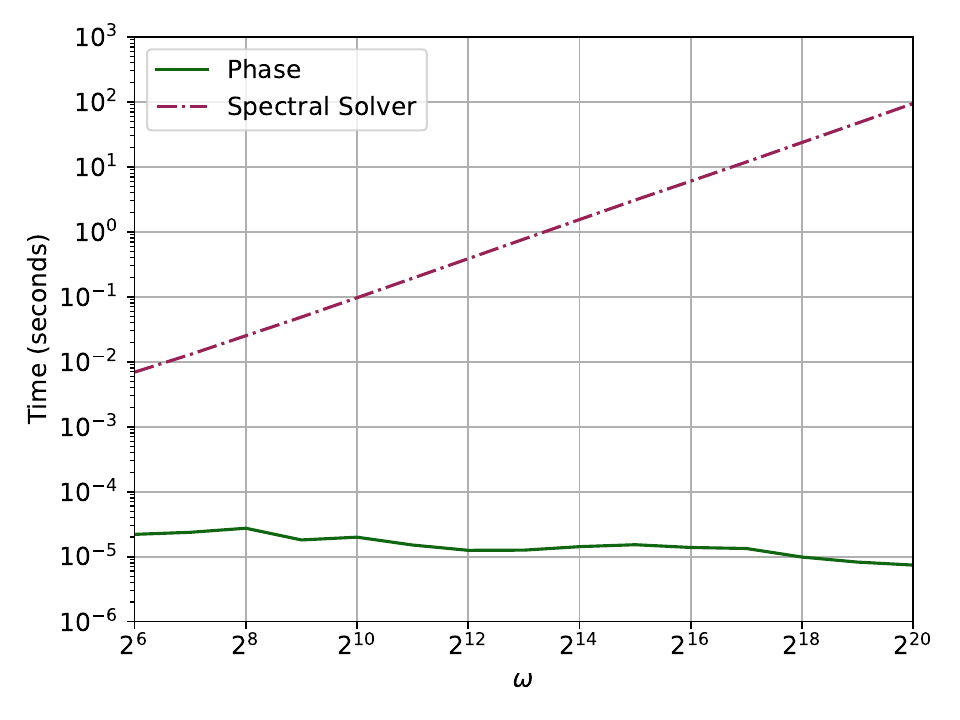}
\hfil
\includegraphics[width=.49\textwidth]{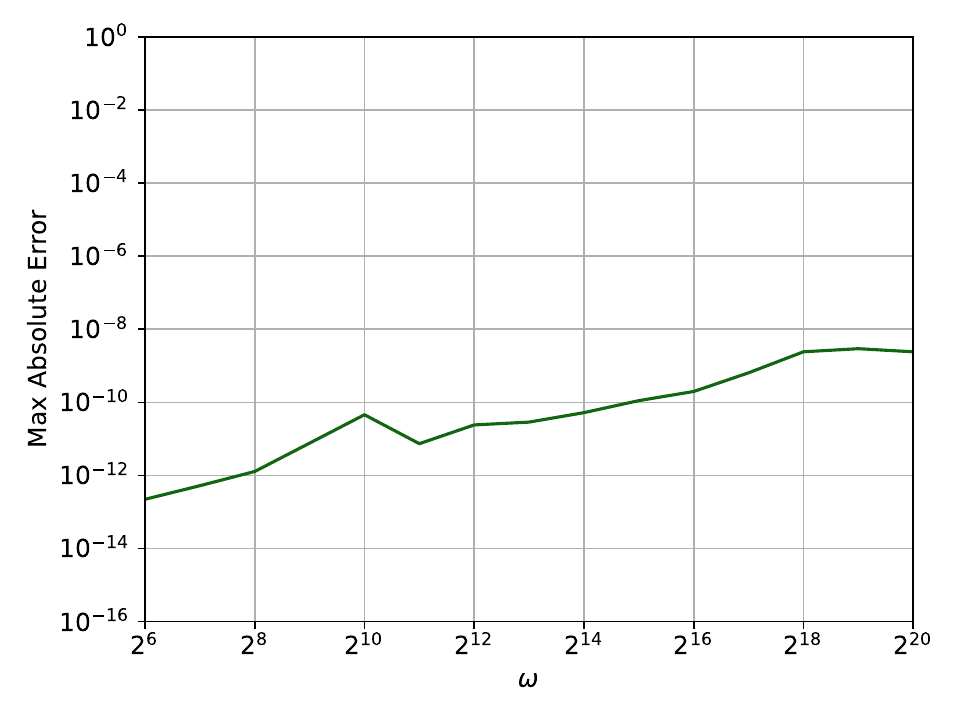}
\hfil

\caption{The results of the experiment of Subsection~\ref{section:experiments5}
in which our algorithm was applied to boundary value problem for
the differential equation (\ref{experiments:bvp})
and the accuracy of the resulting solution was tested using
a standard adaptive Chebyshev spectral solver.
The plot on the left  gives the time required to solve the problem
using each method as a function of the frequency $\omega$,
while the plot on the right gives the maximum 
absolute error observed while comparing the  solution obtained using our method
with the reference solution computed via the standard solver.
}

\label{experiments5:figure1}
\end{figure}

\end{subsection}

%
%


\end{section}

\begin{section}{Conclusion}
\label{section:conclusion}

We have introduced a novel method for solving oscillatory second order linear ordinary
differential equations.  The running time of our scheme is independent of frequency, 
and it is considerably faster than previous methods with this property.  
In the high-frequency regime, it applies a remarkably simple approach, namely, 
discretizing the Riccati equation  via a standard Chebyshev spectral method and 
using the Newton-Kantorovich algorithm to invert the resulting linear system.
Moreover, the results of Section~\ref{section:nonoscillatory}
and (\ref{section:riccati}) rigorously establish the validity of this method.

The machinery necessary to handle the low-frequency regime is somewhat more complicated, 
and this is a weakness of the algorithm of this paper which the authors hope to address 
shortly.      In most cases of interest, the low-frequency regime is close to a turning point
of the differential equation.  Although phase functions can be used to represent solutions of second
order linear ordinary differential equations near a turning point, their derivatives  behave like steep errors 
functions there, which   makes them expensive to approximate via polynomial expansions.  
The  authors are currently developing an algorithm based on an alternate approach to representing solutions
of second order linear ordinary differential equations near turning points.
This should allow for the further acceleration of the method discussed
here, as well as making it simpler and more robust.

We also observe that the bulk of the calculations performed by the algorithm of this paper can be 
parallelized.   There are two problems that must be overcome to do so, but both
admit obvious solutions.  
First, the solution of Appell's equation is carried out sequentially, after the Riccati
equation has been solved over each high-frequency interval.  Rather than doing this, though,
we could construct a basis in the space of solutions of Appell's equation for each
low-frequency interval and construct the desired phase function from these basis
functions after either Appell's equation or the Riccati equation has been solved 
in each interval.  This would allow the bulk of the computations to be 
parallelized.  However, another problem remains, namely, we adaptively
discretize the phase function as we compute it and such calculations are difficult
to parallelize.  It would behoove us to have a set of discretization intervals
determined in advance, before solving the Riccati equation or Appell's equation.
Fortunately, there is a relatively  straightforward mechanism for determining
a suitable set of  discretization intervals {\it a priori}, namely,
we can adaptively discretize the asymptotic approximate  $\rlg$.  The authors
have found that this provides a collection of discretization intervals
which suffices to represent the phase function in almost all cases of interest.

Finally, we note that because many special functions satisfy second order equations,
there are numerous applications of this work to special functions.
For instance, it should allow for the extremely rapid and high-accuracy
computation of classical Gaussian quadrature rules.  
The authors are also interested in using the techniques discussed
here to solve Sturm-Liouville problems, and to rapidly apply the associated
eigentransforms.

\end{section}

\begin{section}{Acknowledgments}
JB  was supported in part by NSERC Discovery grant  RGPIN-2021-02613.
The authors thank Kirill Serkh for several useful discussions.
\end{section}


\bibliographystyle{acm}
\bibliography{riciter.bib}

\begin{thebibliography}{10}

\bibitem{ARDC}
{\sc Agocs, F.~J., and Barnett, A.~H.}
\newblock An adaptive spectral method for oscillatory second-order linear odes
  with frequency-independent cost.
\newblock {\em SIAM Journal on Numerical Analysis 62}, 1 (2024), 295--321.

\bibitem{Agocs}
{\sc Agocs, F.~J., Handley, W.~J., Lasenby, A.~N., and Hobson, M.~P.}
\newblock Efficient method for solving highly oscillatory ordinary differential
  equations with applications to physical systems.
\newblock {\em Phys. Rev. Res. 2\/} (Jan 2020), 013030.

\bibitem{appell}
{\sc Appell, P.}
\newblock Sur la transformation des \'equations diff\'erentielles lin\'eaires.
\newblock {\em Comptes Rendus 91\/} (1880), 211--214.

\bibitem{Arnold}
{\sc Arnold, A., Abdallah, N.~B., and Negulescu, C.}
\newblock {WKB}-based schemes for the oscillatory {1D} {S}chr\"odinger equation
  in the semiclassical limit.
\newblock {\em SIAM Journal on Numerical Analysis 49}, 4 (2011), 1436--1460.

\bibitem{HTFI}
{\sc Bateman, H., and Erd\'elyi, A.}
\newblock {\em Higher Transcendental Functions}, vol.~I.
\newblock McGraw-Hill, New York, New York, 1953.

\bibitem{BremerZeros}
{\sc Bremer, J.}
\newblock On the numerical calculation of the roots of special functions
  satisfying second order ordinary differential equations.
\newblock {\em SIAM Journal on Scientific Computing 39\/} (2017), A55--A82.

\bibitem{BremerPhase}
{\sc Bremer, J.}
\newblock On the numerical solution of second order differential equations in
  the high-frequency regime.
\newblock {\em Applied and Computational Harmonic Analysis 44\/} (2018),
  312--349.

\bibitem{BremerPhase2}
{\sc Bremer, J.}
\newblock Phase function methods for second order linear ordinary differential
  equations with turning points.
\newblock {\em Applied and Computational Harmonic Analysis 65\/} (2023),
  137--169.

\bibitem{BremerSHT}
{\sc Bremer, J., Chen, Z., and Yang, H.}
\newblock Rapid application of the spherical harmonic transform via
  interpolative decomposition butterfly factorization.
\newblock {\em SIAM Journal on Scientific Computing 43}, 6 (2021),
  A3789--A3808.

\bibitem{Ciarlet}
{\sc Ciarlet, P.}
\newblock {\em Linear and Nonlinear Functional Analysis with Applications}.
\newblock Society for Industrial and Applied Mathematics, Philadelphia, PA,
  2013.

\bibitem{Davidson}
{\sc Davidson, R., and Hong, Q.}
\newblock {\em Physics of intense charged particle beams in high energy
  accelerators}.
\newblock World Scientific, Singapore, 2001.

\bibitem{DLMF}
{\it NIST Digital Library of Mathematical Functions}.
\newblock http://dlmf.nist.gov/, Release 1.1.0 of 2020-12-15.
\newblock F.~W.~J. Olver, A.~B. {Olde Daalhuis}, D.~W. Lozier, B.~I. Schneider,
  R.~F. Boisvert, C.~W. Clark, B.~R. Miller, B.~V. Saunders, H.~S. Cohl, and
  M.~A. McClain, eds.

\bibitem{durand75}
{\sc Durand, L.}
\newblock {N}icholson-type integrals for products of {G}egenbauer functions and
  related topics.
\newblock In {\em Theory and Application of Special Functions}, R.~A. Askey,
  Ed. Academic Press, 1975, pp.~353--374.

\bibitem{Einaudi}
{\sc Einaudi, F., and Hines, C.}
\newblock {WKB} approximation in application to acoustic‐gravity waves.
\newblock {\em Canadian Journal of Physics 48\/} (02 2011), 1458--1471.

\bibitem{Hazeltine}
{\sc Hazeltine, R.~D., and Meiss, J.~D.}
\newblock {\em Plasma confinement}.
\newblock Courier Corporation, North Chelmsford, Massachusett, 2003.

\bibitem{BremerRokhlin}
{\sc Heitman, Z., Bremer, J., and Rokhlin, V.}
\newblock On the existence of nonoscillatory phase functions for second order
  ordinary differential equations in the high-frequency regime.
\newblock {\em Journal of Computational Physics 290\/} (2015), 1--27.

\bibitem{Horn}
{\sc Horn, R.~A., and Johnson, C.~R.}
\newblock {\em Matrix Analysis}.
\newblock Cambridge University Press, Cambridge, England, 1990.

\bibitem{ModMagnus}
{\sc Iserles, A.}
\newblock On the global error of discretization methods for highly-oscillatory
  ordinary differential equations.
\newblock {\em BIT Numerical Mathematics 32\/} (2002), 561--599.

\bibitem{Iserles2}
{\sc Iserles, A.}
\newblock Think globally, act locally: solving highly-oscillatory ordinary
  differential equations.
\newblock {\em Applied Numerical Mathematics 43\/} (2002), 145--160.

\bibitem{Kantorovich}
{\sc Kantorovich, L.}
\newblock Functional analysis and applied mathematics.
\newblock {\em Uspehi Matematiceskii Nauk 3\/} (1948), 89--185.

\bibitem{Korner}
{\sc K\"rner, J., Arnold, A., and Döpfner, K.}
\newblock {WKB}-based scheme with adaptive step size control for the
  {S}chr\"odinger equation in the highly oscillatory regime.
\newblock {\em Journal of Computational and Applied Mathematics 404\/} (2022),
  113905.

\bibitem{Kummer}
{\sc Kummer, E.}
\newblock De generali quadam aequatione differentiali tertti ordinis.
\newblock {\em Progr. Evang. K\"ongil. Stadtgymnasium Liegnitz\/} (1834).

\bibitem{Lubich}
{\sc Lorenz, K., Jahnke, T., and Lubich, C.}
\newblock Adiabatic integrators for highly oscillatory second-order linear
  differential equations with time-varying eigendecomposition.
\newblock {\em BIT Numerical Mathematics\/} (2005), 91--115.

\bibitem{Jerome}
{\sc Martin, J., and Schwarz, D.}
\newblock {WKB} approximation for inflationary cosmological perturbations.
\newblock {\em Physical Review D 67\/} (10 2002).

\bibitem{Olver}
{\sc Olver, F.~W.}
\newblock {\em Asymptotics and Special Functions}.
\newblock A.K. Peters, Wellesley, Massachusetts, 1997.

\bibitem{SheehanOlver2}
{\sc Olver, S.}
\newblock {GMRES} for the differentiation operator.
\newblock {\em SIAM Journal on Numerical Analysis 47}, 5 (2009), 3359--3373.

\bibitem{SheehanOlver1}
{\sc Olver, S.}
\newblock {GMRES} for oscillatory matrix-valued differential equations.
\newblock In {\em Spectral and High Order Methods for Partial Differential
  Equations\/} (Berlin, Heidelberg, 2011), J.~S. Hesthaven and E.~M.
  R{\o}nquist, Eds., Springer Berlin Heidelberg, pp.~267--274.

\bibitem{2018_Pritula}
{\sc Pritula, G.~M., Petrenko, E.~V., and Usatenko, O.~V.}
\newblock Adiabatic dynamics of one-dimensional classical {H}amiltonian
  dissipative systems.
\newblock {\em Physics Letters, Section A: General, Atomic and Solid State
  Physics 382}, 8 (feb 2018), 548--553.

\bibitem{SpiglerPhase2}
{\sc Spigler, R.}
\newblock Asymptotic-numerical approximations for highly oscillatory
  second-order differential equations by the phase function method.
\newblock {\em Journal of Mathematical Analysis and Applications 463\/} (2018),
  318--344.

\bibitem{SpiglerZeros}
{\sc Spigler, R., and Vianello, M.}
\newblock A numerical method for evaluating the zeros of solutions of
  second-order linear differential equations.
\newblock {\em Mathematics of Computation 55\/} (1990), 591--612.

\bibitem{SpiglerPhase1}
{\sc Spigler, R., and Vianello, M.}
\newblock The phase function method to solve second-order asymptotically
  polynomial differential equations.
\newblock {\em Numerische Mathematik 121\/} (2012), 565--586.

\bibitem{trefethen}
{\sc Trefethen, L.}
\newblock Is {G}auss quadrature better than {C}lenshaw–{C}urtis?
\newblock {\em SIAM Review 50\/} (2008), 67--87.

\end{thebibliography}

\end{document}